\documentclass[a4paper, 12pt,twoside]{article}
\usepackage[english]{babel}
\usepackage[latin1]{inputenc}
\usepackage{amsmath}
\usepackage{amssymb,amsfonts}
\usepackage{graphicx}
\usepackage{times,amssymb,amscd}
\newtheorem{thm}{Theorem}[section]

\newtheorem{rem}[thm]{Remark}
\numberwithin{equation}{section}

\newcommand{\bA}{\mathbf{A}}

\newcommand{\bE}{\mathbf{E}}
\newcommand{\bG}{\mathbf{G}}
\newcommand{\bH}{\mathbf{H}}

\newcommand{\bL}{\mathbf{L}}

\newcommand{\bR}{\mathbf{R}}

\newcommand{\bS}{\mathbf{S}}
\newcommand{\bV}{\mathbf{V}}
\newcommand{\bZ}{\mathbf{Z}}

\newcommand{\be}{\mathbf{e}}

\newcommand{\bx}{\mathbf{x}}
\newcommand{\by}{\mathbf{y}}

\newcommand{\bg}{\mathbf{g}}

\newcommand{\bI}{\mathbf{I}}

\newcommand{\BV}{\boldsymbol{V}}

\newcommand{\Be}{\boldsymbol{e}}
\newcommand{\Bu}{\boldsymbol{u}}
\newcommand{\Bv}{\boldsymbol{v}}

\newcommand{\cD}{\mathcal{D}}

\newcommand{\cP}{\mathcal{P}}
\newcommand{\cS}{\mathcal{S}}

\newcommand{\cB}{\mathcal{B}}

\newcommand{\HYP}{\bH^3}
\newcommand{\SXR}{\bS^2\!\times\!\bR}
\newcommand{\HXR}{\bH^2\!\times\!\bR}
\newcommand{\SLR}{\widetilde{\bS\bL_2\bR}}
\newcommand{\NIL}{\mathbf{Nil}}
\newcommand{\SOL}{\mathbf{Sol}}

\begin{document}
\pagestyle{myheadings}
\markboth{\centerline{Arnasli Yahya and Jen\H o Szirmai }}
{New lower bound...}
\title
{New lower bound for the optimal congruent geodesic ball packing density of screw motion groups in $\mathbf{H}^2\!\times\!\mathbf{R}$ space. 
}

\author{Arnasli Yahya and Jen\H o Szirmai \\
\normalsize Department of Algebra and Geometry, Institute of Mathematics,\\
\normalsize Budapest University of Technology and Economics, \\
\normalsize M\H uegyetem rkp. 3., H-1111 Budapest, Hungary \\
\normalsize arnasli@math.bme.hu,~szirmai@math.bme.hu
\date{\normalsize{\today}}}
\maketitle


\maketitle
\begin{abstract}
In this paper, we present a new record for the densest geodesic congruent ball packing configurations in $\HXR$ geometry, generated by screw motion groups. 
These groups are derived from the direct product of rotational groups on  
$\mathbf{H}^2$ and some translation components on the real fibre direction $\mathbf{R}$ that can be determined by the corresponding Frobenius congruences. 
Moreover, we developed a procedure to determine the optimal radius for the densest geodesic ball packing configurations related to the considered screw motion groups.

The highest packing density, $\approx0.80529$, is achieved by a multi-transitive case given by rotational parameters $(2,20,4)$. 

E. Moln\'{a}r demonstrated that homogeneous 3-spaces can be uniformly interpreted in the projective 3-sphere 
$\mathcal{PS}^3(\bV^4, \BV_4, \mathbf{R})$. We use this projective model of $\HXR$ to compute and visualize the locally optimal geodesic ball arrangements.
\end{abstract}

\newtheorem{Theorem}{Theorem}[section]
\newtheorem{corollary}[Theorem]{Corollary}
\newtheorem{lemma}[Theorem]{Lemma}
\newtheorem{exmple}[Theorem]{Example}
\newtheorem{definition}[Theorem]{Definition}
\newtheorem{Remark}[Theorem]{Remark}
\newtheorem{proposition}[Theorem]{Proposition}
\newenvironment{remark}{\begin{rmrk}\normalfont}{\end{rmrk}}
\newenvironment{example}{\begin{exmple}\normalfont}{\end{exmple}}
\newenvironment{acknowledgement}{Acknowledgement}

\section{Introduction}\label{sec_1}
The second author extended the classic Kepler problem to non-constant curvature 
Thurston geometries $\SXR,~\HXR,~$ $\SLR,~\NIL,~\SOL$, 
in \cite{T}, \cite{S}, \cite{MSz}, \cite{Sz13-1}. The investigation of this issue brought many interesting results and opened an important 
path in the direction of non-Euclidean crystal geometry (see the survey \cite{Sz22-3} 
and \cite{Sz11-1}, \cite{Sz11-2}, \cite{MSz14},\cite{MSzV17},\cite{Sz23}, \cite{FTL}, \cite{G--K--K}). 
We mention only some here:
\begin{enumerate}
\item In \cite{Sz07} we investigated the geodesic balls of the $\NIL$ space and computed their volume,
introduced the notion of the $\NIL$ lattice, $\NIL$ parallelepiped and the density of the lattice-like ball packing.
Moreover, we determined the densest lattice-like geodesic ball packing. The density of this densest packing is
$\approx 0.7809$, which may be quite surprising
when compared with the Euclidean result $\frac{\pi}{\sqrt{18}} \approx 0.74048$. The kissing number of the balls
in this packing is $14$.
\item Moreover, a candidate of the densest geodesic ball
packing is described in \cite{Sz13-1}. In the Thurston geometries the greatest known density was $\approx 0.8533$ which is not realized by a packing with {\it equal balls} of the hyperbolic
space $\HYP$. However, it is attained by, e.g., by a {\it horoball packing} of
$\overline{\bH}^3$ where the ideal centres of horoballs lie on the
absolute figure of $\overline{\bH}^3$ inducing the regular ideal
simplex tiling $(3,3,6)$ by its Coxeter-Schl\"afli symbol.
In \cite{Sz13-1} we have presented a geodesic ball packing in the $\SXR$ geometry
whose density is $\approx 0.8776$.
\item  In \cite{Sz12-5} we determined the geodesic balls of $\HXR$ and computed their volume, and
 defined the notion of the geodesic ball packing and its density.
 Moreover, we have developed a procedure to determine the density of the simply or multiply transitive geodesic ball packings for
 generalized Coxeter space groups of $\HXR$ and applied this algorithm to them.
 The Dirichlet-- Voronoi cells for the above space groups are ``prisms" in the $\HXR$ sense.
 The optimal packing density of the generalized Coxeter space groups is $\approx 0.60726$. 
\item In the paper \cite{YSz23-1}, we investigated the locally optimal geodesic congruent ball packings generated by rotational point groups where the corresponding translation parts 
are trivial. Moreover, we studied the monotonicity behavior of their densities 
in $\mathbf{H}^2\!\times\!\mathbf{R}$ space. 
 
\end{enumerate}
Our article is related to the previous work, in which we study the locally optimal geodesic ball packings 
with equal balls to the space groups having rotation point groups and generators that are proper screw motions in $\HXR$ geometry. 

The findings of this study are summarized in Theorems 3.1, 3.3, 3.9, and 3.10. The numerical results are systematically presented in Tables 1-8. Notably, we have determined that the optimal packing density, $\approx 0.80529$, exceeds the densities reported in previous studies. Consequently, this result establishes a new thinnest lower bound for the optimal density of ball packing. We posit that the constructive methods outlined in this paper can be extended to explore geodesic ball packings for general space groups in $\HXR$.\\ In the following section, we begin by reviewing the structure of the space group, the geodesic curves, and the balls in $\HXR$.
\section{On $\HXR$ geometry}
\subsection{The structure of $\HXR$ space groups}
$\HXR$ is one of the eight simply connected 3-dimensional maximal homogeneous Riemannian geometries. 
This Seifert fibre space is derived by the direct product of the hyperbolic plane $\bH^2$ and the real line $\bR$. 
The points are described by $(P,p)$ where $P\in \bH^2$ and $p\in \bR$. The complete isometry group $Isom(\HXR)$ of $\HXR$ can be 
derived by the direct product of the isometry group $Isom(\bH^2)$ of the hyperbolic plane and the isometry group $Isom(\bR)$ of the real line 
as follows (see \cite{F01}, \cite{Sz12-5}). 
\begin{equation}
\begin{gathered}
Isom(\HXR):=Isom(\bH^2) \times Isom(\bR);\\
Isom(\bH^2):=\{A  \ : \ \bH^2 \mapsto \bH^2 \ : \ (P,p) \mapsto (PA,p) \} \ \text{for any fixed $p \in \bR$}.  \\
Isom(\bR):=\{\rho ~ : ~ (P,p) \mapsto (P, \pm p + \tau) \}, \ \text{for any fixed $P\in \bH^2$}.\\
\text{here the "-" sign provides a reflection in the point} \ \frac{\tau}{2} \in \bR, \\ \text{by the "+" sign we get a translation of $\bR$}.  
\end{gathered} \tag{2.1}
\end{equation}

The structure of a discontinuously acting, so finitely generated isometry group $\Gamma \subset$ $Isom(\HXR)$ is the following (see \cite{Sz12-5}:
$\Gamma:=\langle (A_1 \times \rho_1), \dots (A_n \times \rho_n) \rangle$, where
$A_i \times \rho_i:=A_i \times (R_i,\tau_i):=(g_i,\tau_i)$, $(i \in \{ 1,2, \dots n \}$ and $A_i \in Isom(\bH^2)$, $R_i$ is either the identity map 
$\mathbf{1_R}$ of $\bR$ or the point reflection $\overline{\mathbf{1}}_{\mathbf{R}}$. $g_i:=A_i \times R_i$ is called the linear part of the transformation
$(A_i \times \rho_i)$ and $\tau_i$ is its translation part. 
The multiplication formula is the following:
\begin{equation}
(A_1 \times R_1,\tau_1) \circ (A_2 \times R_2,\tau_2)=((A_1A_2 \times R_1R_2,\tau_1R_2+\tau_2). \tag{2.2}
\end{equation}
\begin{definition}
$L_{\Gamma}$ is a \textbf{one dimensional lattice} on $\bR$ fibres if there is a positive real number $r$ such that     
\begin{equation}
L_{\Gamma}:= \{ kr:(P,p) \mapsto (P,p+kr), \ \forall P \in \mathbf{H}^2; \ \forall p \in \bR \ | \  0 < r \in \bR, \ k \in \bZ \} \notag
\end{equation}
\end{definition}
\begin{definition}
A group of isometries $\Gamma \subset Isom(\HXR)$ is called a \textbf{space group} if its linear parts form a cocompact 
(i.e. of compact fundamental domain in $\bH^2$) group $\Gamma_0$ called the point group of 
$\Gamma$, moreover, the translation parts to the identity of this point group are required to form a one dimensional lattice $L_{\Gamma}$ of $\bR$.
\end{definition}
\begin{Remark}
\begin{enumerate}
\item It can easily be proved, that such a space group $\Gamma$ has a compact fundamental domain $\mathcal{F}_\Gamma$ in $\HXR$.  
\end{enumerate}
\end{Remark}
\begin{definition}
The $\HXR$ space groups $\Gamma_1$ and $\Gamma_2$ are geometrically equivalent, called \textbf{equivariant}, if there is a "similarity" transformation
$\Sigma:=S \times \sigma$ $(S \in {Hom}(\bH^2), \sigma \in Sim(\bR))$, such that $\Gamma_2=\Sigma^{-1} \Gamma_1 \Sigma$, where 
$S$ is a piecewise linear (i.e. $PL$) homeomorphism of $\bH^2$ which deforms the fundamental domain of $\Gamma_1$ into that of $\Gamma_2$.
Here $\sigma(s,t):p \rightarrow p \cdot s+t$ is a similarity of $\bR$, i.e. multiplication by $0 \ne s \in \bR$ and then addition by $t \in \bR$ for every 
$p \in \bR$. 
\end{definition}

The equivariance class of a hyperbolic plane group or its orbifold \cite{LM91} can be characterized by its {\it Macbeath-signature}. 
In 1967-69 {\sc Macbeath} completed the classification of hyperbolic crystallographic plane groups, (for short NEC groups) \cite{M}. 
He considered isometries containing orientation -preserving and -reversing transformations as well in the hyperbolic plane. His paper deals
with  NEC groups, but the Macbeath-signature economically characterizes the Euclidean and spherical plane groups, too. The signature of 
a plane group is the following
\begin{equation}
(\pm,g;[m_1,m_2, \dots, m_r];~\{(n_{11},n_{12},\dots,n_{1s_1}), \dots , (n_{k1},n_{k2},\dots,n_{ks_k})\}). \tag{2.3}
\end{equation}
and, with the same notations, the combinatorial measure $T$ of the fundamental polygon is expressed by:
\begin{equation}
T\kappa=\pi\Big\{\sum_{l=1}^{r} \Big(\frac{2}{m_l}-2 \Big)+\sum_{i=1}^{k}\Big(-2+\sum_{j=s_1}^{s_i}\Big(-1+\frac{1}{n_{ij}}\Big)\Big)+2\chi\Big\}. \notag  
\end{equation}
Here $\chi=2-\alpha g$ $(\alpha=1$ for $-$, $\alpha=2$ for $+$, the sign $\pm$ refers to orientability) $\chi$ is the Euler characteristic of the surface
with genus $g$, and $\kappa$ will denote the Gaussian curvature of the realizing plane $\bS^2$, $\bE^2$ 
or $\bH^2$, whenewer $\kappa>0$, $\kappa=0$ or $\kappa<0$, respectively. The genus $g$, the proper periods $m_l$ of $r$ rotation centres and the period-cycles
$(n_{i1},n_{i2},\dots,n_{is_i})$ of dihedral corners on $i^{th}$ one of the $k$ boundary components, together, 
with a marked fundamental polygon with side pairing generators and with a corresponding group presentation determine a plane group up 
to a well formulated equivariance for $\bS^2$, $\bE^2$ 
and $\bH^2$, respectively \cite{M}. 
\begin{Theorem}[\cite{Sz12-5}]
Let $\Gamma$ be a $\HXR$ space group, its point group $\Gamma_0$ belongs to one of the following three types:
\begin{enumerate}
\item[I.] $\bG_{\bf H^2} \times {\bf 1}_{\bf R}$, 
${\bf 1}_{\bf R} : x \mapsto x$ is the identity of $\bR$.
\item[II.] $\bG_{\bf H^2} \times \langle \overline{{\bf 1}}_{\bf R} \rangle$, where $\overline{\bf 1}_{\bR}: x \mapsto -x+r$ is the $\frac{r}{2}$ reflection
of $\bR$ with some $r$ and $\langle \overline{{\bf 1}}_{\bR} \rangle$ denotes its special linear group of two elements. 
\item[III.] If the hyperbolic group $\bG_{\bf H^2}$ contains a normal subgroup $\bG$ of index two, then 
$\bG_{\bf H^2}\bG:=\{\bG \times {\bf 1}_{\bf R}\} \cup \{(\bG_{\bf H^2} \setminus \bG) \times \overline{{\bf 1}}_\bR\}$ forms a point group.
\end{enumerate}
Here $\bG_{\bf H^2}$ is a group of hyperbolic isometries with compact fundamental domain $\mathcal{F}_\Gamma$. 
\end{Theorem}
In this paper we consider space groups having rotation point groups
and their generators are screw motions in $\HXR$ geometry 
\begin{definition}
A $\HXR$ space group $\Gamma$ is called {\textbf{generalized screw motion group}} if the generators $\bg_i, \ (i=1,2,\dots m)$ of its point group $\Gamma_0$ 
are rotations and the possible translation parts of all the above generators are lattice translations, i.e. $\tau_i \ \text{mod}L_\Gamma \ (i=1,2).$
\end{definition}

In this paper we deal with ``generalized screw motion groups" in $\HXR$ space given by parameters $ 3 \le p_1,p_2 \in \mathbb{N}$ where 
$\frac{1}{p_1}+\frac{1}{p_2} < \frac{1}{2}$,
\begin{equation}
\Gamma_{(p_0=2,p_1,p_2)} ~ (+,~0,~ [p_0=2,p_1,p_2];~\{~ \}), \ \ \Gamma_0=(\bg_0,\bg_1 - \bg_0^{2},\bg_1^{p_1}, (\bg_{0}\bg_1)^{p_2}). \tag{2.4} 
\end{equation}
For a fundamental domain of the above space groups we can combine a fundamental domain of a rotation group of the hyperbolic plan 
with a part of a real line segment $r$. 
\subsection{Geodesic curves and balls in $\HXR$}\label{sect3}
In \cite{M97} E. {Moln\'ar} has shown, that the homogeneous 3-spaces
have a unified interpretation in the projective 3-sphere $\mathcal{PS}^3(\bV^4,\BV_4, \mathbf{R})$. 
In our work we shall use this projective model of $\HXR$ and
the Cartesian homogeneous coordinate simplex $E_0(\be_0)$, $E_1^{\infty}(\be_1)$, $E_2^{\infty}(\be_2)$,
$E_3^{\infty}(\be_3)$, $(\{\be_i\}\subset \bV^4$)$ \text{ with the unit point}$ $E(\be = \be_0 + \be_1 + \be_2 + \be_3 ))$ 
which is distinguished by an origin $E_0$ and by the ideal points of coordinate axes, respectively. 
Moreover, $\by=c\bx$ with $0<c\in \mathbf{R}$ (or $c\in\mathbb{R}\setminus\{0\})$
defines a point $(\bx)=(\by)$ of the projective 3-sphere $\cP \cS^3$ (or that of the projective space $\cP^3$ where opposite rays
$(\bx)$ and $(-\bx)$ are identified). 
The dual system $\{(\Be^i)\}\subset \BV_4$ describes the simplex planes, especially the plane at infinity 
$(\Be^0)=E_1^{\infty}E_2^{\infty}E_3^{\infty}$, and generally, $\Bv=\Bu\frac{1}{c}$ defines a plane $(\Bu)=(\Bv)$ of $\cP \cS^3$
(or that of $\cP^3$). Thus $0=\bx\Bu=\by\Bv$ defines the incidence of point $(\bx)=(\by)$ and plane
$(\Bu)=(\Bv)$, as $(\bx) \text{I} (\Bu)$ also denotes it. Thus {$\HXR$} can be visualized in the affine 3-space $\bA^3$
(so in $\bE^3$) as well.

The points of $\HXR$ space, forming an open cone solid in the projective space $\mathcal{P}^3$, are the following:
\begin{equation}
\HXR:=\big\{ X(\bx=x^i \be_i)\in \mathcal{P}^3: -(x^1)^2+(x^2)^2+(x^3)^2<0<x^0,~x^1 \big\}. \notag
\end{equation}
In this context E. Moln\'ar \cite{M97} has derived the infinitezimal arc-length square at any point of $\HXR$ as follows
\begin{equation}
   \begin{gathered}
     (ds)^2=\frac{1}{(-x^2+y^2+z^2)^2}\cdot \Big\{[(x)^2+(y)^2+(z)^2](dx)^2+ \\ + 2dxdy(-2xy)+2dxdz (-2xz)+ [(x)^2+(y)^2-(z)^2] (dy)^2+ \\ 
     +2dydz(2yz)+ [(x)^2-(y)^2+(z)^2](dz)^2\Big\}.
       \end{gathered} \tag{2.5}
     \end{equation}
This becomes simpler in the following special (cylindrical) coordinates $(t, r, \alpha), \ \ (r \ge 0, ~ -\pi < \alpha \le \pi)$ 
with the fibre coordinate $t \in \bR$. We describe points in our model by the following equations: 
\begin{equation}
x^0=1, \ \ x^1=e^t \cosh{r},  \ \ x^2=e^t \sinh{r} \cos{\alpha},  \ \ x^3=e^t \sinh{r} \sin{\alpha}  \tag{2.6}.
\end{equation}
Then we have $x=\frac{x^1}{x^0}=x^1$, $y=\frac{x^2}{x^0}=x^2$, $z=\frac{x^3}{x^0}=x^3$, i.e. the usual Cartesian coordinates.
We obtain by \cite{M97} that in this parametrization the infinitezimal arc-length square and the symmetric metric tensor field $g_{ij}$ by (3.1):
at any point of $\HXR$ is the following
     \begin{equation}
       g_{ij}:=
       \begin{pmatrix}
         1&0&0 \\
         0&1 &0 \\
         0&0& \sinh^2{r} \\
         \end{pmatrix}. \tag{2.7}
     \end{equation}
The geodesic curves of $\HXR$ are generally defined as having locally minimal arc length between their two arbitrary (near enough) points. 
The equation systems of the parametrized geodesic curves $\gamma(t(\tau),r(\tau),\alpha(\tau))$ in our model can be determined by the 
general theory of Riemannian geometry:

We can assume, that the starting point of a geodesic curve is $(1,1,0,0)$, because we can transform a curve into an 
arbitrary starting point, moreover, unit velocity with "geographic" coordinates $(u,v)$ can be assumed:
\begin{equation}
\begin{gathered}
        r(0)=\alpha(0)=t(0)=0; \ \ \dot{t}(0)= \sin{v}, \ \dot{r}(0)=\cos{v} \cos{u}, \dot{\alpha}(0)=\cos{v} \sin{u}; \\
        - \pi < u \leq \pi, ~ -\frac{\pi}{2}\le v \le \frac{\pi}{2}. \notag
\end{gathered}
\end{equation}
Then by (3.2) we get with $c=\sin{v}$, $\omega=\cos{v}$ the equation systems of a geodesic curve:
\begin{equation}\label{eq:geodesic}
  \begin{gathered}
   x(\tau)=e^{\tau \sin{v}} \cosh{(\tau \cos{v})}, \\ 
   y(\tau)=e^{\tau \sin{v}} \sinh{(\tau \cos{v})} \cos{u}, \\
   z(\tau)=e^{\tau \sin{v}} \sinh{(\tau \cos{v})} \sin{u},\\
   -\pi < u \le \pi,\ \ -\frac{\pi}{2}\le v \le \frac{\pi}{2}. \tag{2.8}
  \end{gathered}
\end{equation}
\begin{definition}
The \textbf{distance} $d(P_1,P_2)$ between the points $P_1$ and $P_2$ is defined by the arc length of the geodesic curve 
from $P_1$ to $P_2$.
\end{definition}
 \begin{definition}
 The \textbf{geodesic sphere} of radius $\rho$ (denoted by $S_{P_1}(\rho)$) with centre at the point $P_1$ is defined as the set of all points 
 $P_2$ in the space with the condition $d(P_1,P_2)=\rho$. Moreover, we require that the geodesic sphere is a simply connected 
 surface without selfintersection in $\HXR$ space.
 \end{definition}
 \begin{Remark}
 In this paper we consider only the usual spheres with "proper centre", i.e. $P_1 \in \HXR$. 
 If the centre of a "sphere" lie on the absolute quadric or lie out of our model the notion of the "sphere" (similarly to the hyperbolic space),
 can be defined, but that cases we shall study in a forthcoming work.
 \end{Remark}
 \begin{definition}
 The body of the geodesic sphere of centre $P_1$ and of radius $\rho$ in $\HXR$ space is called \textbf{geodesic ball}, denoted by $B_{P_1}(\rho)$,
 i.e. $Q \in B_{P_1}(\rho)$ iff $0 \leq d(P_1,Q) \leq \rho$.
 \end{definition}
In \cite{Sz12-5} we determined the volume of a geodesic ball: 

\begin{equation}\label{Eq:volume geodesic ball}
\begin{gathered}
Vol(B(\rho))=\int_{V} \frac{1}{(x^2-y^2-z^2)^{3/2}}\mathrm{d}x ~ \mathrm{d}y ~ \mathrm{d}z = \\ = \int_{0}^{\rho} \int_{-\frac{\pi}{2}}^{\frac{\pi}{2}} 
\int_{-\pi}^{\pi} 
|\tau \cdot \sinh(\tau \cos(v))| ~ \mathrm{d} u \ \mathrm{d}v \ \mathrm{d}\tau = \\ =
2 \pi \int_{0}^{\rho} \int_{-\frac{\pi}{2}}^{\frac{\pi}{2}} |\tau \cdot \sinh(\tau \cos(v))| ~ \mathrm{d} v \ \mathrm{d}\tau. \tag{2.9}
\end{gathered}
\end{equation}
We will use for the computations the notion and the volumes of $\HXR$ prisms. The $\HXR$ prism (see \cite{Sz12-5}) is a convex hull of two congruent
$p$-gons $(p > 2)$ in ``parallel planes", (a "plane" is one sheet of concentric two sheeted hyperboloids in our model) related by translation along 
the radii joining their corresponding vertices that are the common perpendicular lines of the two "hyperboloid-planes". 
The prism is a polyhedron having at each vertex one hyperbolic $p$-gon and two "quadrangles".  
The $p$-gonal faces of a prism called cover-faces, and the other faces 
are the side-faces. In these cases every face of each polyhedron meets only one face of another polyhedron. 

The volume of a $\HXR$ $p$-gonal prism can be computed by the following formula:
\begin{equation}\label{Eq: Volume of Cell}
Vol(\mathcal{P})=\mathcal{A} \cdot h \tag{2.10}
\end{equation}
where $\mathcal{A}$ is the area of the hyperbolic $p$-gon in base plane and
$h$ is the height of the prism.

\section{Screw motion groups and their generated packing configurations}
A $\HXR$ space group $\Gamma$ has a compact fundamental domain. Typically, the shape of the fundamental domain of a group of $\bH^2$ is not
determined uniquely but the area of the domain is finite and unique
by its combinatorial measure, (see also \cite{YA2023}). Thus the shape of the fundamental domain of a crystallographic group of $\HXR$ is not unique, as well.

In the following let $\Gamma$ be a fixed by {\it screw motions generated} space group of $\HXR$. We
will denote by $d(X,Y)$ the distance of two points $X$, $Y$ by Definition 2.7.
\begin{definition}
We say that the point set
$$
\cD(K)=\{X\in\HXR\,:\,d(K,X)\leq d(K^\bg,X)\text{ for all }\bg\in\ \Gamma\}
$$
is the \textbf{{Dirichlet--Voronoi cell} (D-V~cell)} to $\Gamma$ around the kernel
point $K\in\HXR$.
\end{definition}
\begin{definition}
We say that
$$
\Gamma_X=\{\bg\in\Gamma\,:\,X^\bg=X\}
$$
is the \textbf{stabilizer subgroup} of $X\in\HXR$ in $\Gamma$.
\end{definition}
\begin{definition}
Assume that the stabilizer $\Gamma_K=\bI$ i.e. $\Gamma$ acts simply transitively on
the orbit of a point $K$. Then let $\cB_K$ denote the \textbf{greatest ball}
of centre $K$ inside the D-V cell $\cD(K)$, moreover let $\rho(K)$ denote the
\textbf{radius} of $\cB_K$. It is easy to see that
$$
\rho(K)=\min_{\bg\in\Gamma\setminus\bI}\frac12 d(K,K^\bg).
$$
\end{definition}
The $\Gamma$-images of $\cB_K$ form a \textbf{\emph{ball packing}} $\cB^\Gamma_K$ with centre
points $K^{\Gamma}$.
\begin{definition}\label{def:dens}
The \textbf{{density}} of ball packing $\cB^\Gamma_K$ is
$$
\delta(K)=\frac{Vol(\cB_K)}{Vol(\cD_K)}.
$$
\end{definition}
It is clear that the orbit $K^\Gamma$ and the ball packing $\cB^\Gamma_K$ have the
same symmetry group, moreover this group contains the starting
crystallographic group $\Gamma$:
$$
Sym K^\Gamma=Sym\cB^\Gamma_K\geq\Gamma.
$$
\begin{definition}
We say that the orbit $K^\Gamma$ and the ball packing $\cB^\Gamma_K$ is
\textbf{characteristic} if $Sym K^\Gamma=\Gamma$, else the orbit is not
characteristic.
\end{definition}
\subsection{Simply transitive case}
\emph{Our problem is} to find a
point $K \in\ \HXR$ and the orbit $K^\Gamma$ for $\Gamma$ such that $\Gamma_K=\bI$
and the density $\delta(K)$ of the corresponding ball packing
$\cB^\Gamma(K)$ is maximal. In this case the ball packing $\cB^\Gamma(K)$ is
said to be \emph{optimal.}

The lattice of $\Gamma$ has a free parameter $p(\Gamma)$. Then we have to find the densest ball packing on $K$ for fixed
$p(\Gamma)$, and vary $p$ to get the optimal ball packing.
\begin{equation}\label{Eq:Global_maximum}
\delta^{opt}(\Gamma)=\max_{K, \ p(\Gamma)}(\delta(K)) \tag{3.1}
\end{equation}
Let $\Gamma$ be a fixed by {\it screw motions generated} group.
The stabilizer of $K$ is trivial i.e. we are looking for the optimal kernel point
in a 3-dimensional region, inside of a fundamental domain of $\Gamma$ with free fibre parameter $p(\Gamma)$.
It can be assumed by the homogeneity of $\HXR$, that the fibre coordinate of the center of the optimal ball is zero.
\subsubsection{Possible translation parts: Frobenius congruences}
In this paper, we deal with ``generalized screw motion groups" in $\HXR$ space given by parameters $ 3 \le p_1,p_2 \in \mathbb{N}$ where 
$\frac{1}{p_1}+\frac{1}{p_2} < \frac{1}{2}$ (see (2.4)). In the folowing we denote the group $\Gamma_{(p_0=2,p_1,p_2)}$ with its corresponding parameters $(2,p_1,p_2)$.

For each $(2,p_1,p_2)$, there exist $(\tau_0,\tau_1,\tau_2)$ that satisfy the Frobenius congruences
\begin{equation*}
    (g_0,\tau_0)^2=(g_1,\tau_1)^{p_1}=(g_1 g_0,\tau_1+\tau_0)^{p_2}=(g_2,\tau_2)^{p_2}=\textbf{I}.
\end{equation*}
Then based on those relations, the fibre translation part from kernel point $\boldsymbol{K}:=(K,t)$, where $K \in \mathbf{H}^2$ and $t \in \mathbf{R}$ would satisfy 
\begin{align*}
    &t+2\tau_0 \equiv t\mod \xi,~~~~~~~
    &t+p_1 \tau_1 \equiv t \mod \xi,\\
    &t+p_2 (\tau_1+\tau_0) \equiv t \mod \xi,~~~~~~~
    &t+p_2 \tau_2 \equiv t \mod \xi,
\end{align*}
for some $0 \leq \xi \in \mathbf{R}$. If we write $r_0:=\frac{\tau_0}{\xi}$, $r_1:=\frac{\tau_1}{\xi}$, and $r_2:=\frac{\tau_2}{\xi}$, then we obtain 
\begin{equation}
    2r_0=k_0,~~p_1 r_1=k_1,~~p_2r_2=k_2,~~p_2(r_1+r_0)=k,
\end{equation}
for some integer $k_0, k_1, k_2, k$.\\
We provide all possible translations in the following
\begin{thm}[Frobenius congruence for simply transitive case]\label{Thm.3.1}
For each $(2,p_1,p_2)$ the solutions to Frobenius congruence is satisfied by
    \begin{enumerate}
        \item If both $p_1$ and $p_2$ are \textbf{even} then\\ $(\tau_0,\tau_1,\tau_2)  \sim  \left(0,\frac{m}{gcd(p_1,p_2)},\frac{m}{gcd(p_1,p_2)}\right),\left( \pm \frac{1}{2},0,\pm\frac{1}{2} \right), \left( \pm \frac{1}{2},\pm\frac{1}{2},0 \right) $
        \item If $p_1$ is \textbf{odd} and $p_2$ is \textbf{even} then\\ $(\tau_0,\tau_1,\tau_2)  \sim  \left(0,\frac{m}{gcd(p_1,p_2)},\frac{m}{gcd(p_1,p_2)}\right),\left(\pm \frac{1}{2},0,\pm\frac{1}{2} \right)$
        \item If $p_1$ is \textbf{even} and $p_2$ is \textbf{odd} then\\ $(\tau_0,\tau_1,\tau_2)  \sim  \left(0,\frac{m}{gcd(p_1,p_2)},\frac{m}{gcd(p_1,p_2)}\right), \left( \pm \frac{1}{2},\pm\frac{1}{2},0 \right) $
        \item If both $p_1$ and $p_2$ are \textbf{odd} then\\ $(\tau_0,\tau_1,\tau_2)  \sim  \left(0,\frac{m}{gcd(p_1,p_2)},\frac{m}{gcd(p_1,p_2)}\right)$,
    \end{enumerate}
    where $-\lfloor \frac{gcd(p_1,p_2)}{2} \rfloor \leq m \leq \lfloor \frac{gcd(p_1,p_2)}{2} \rfloor,~m \in \mathbb{Z}$, and $gcd(p_1,p_2)$ is the greatest common divisor of $p_1$ and $p_2$.
\end{thm}
\textbf{Proof:}\\
We provide the proof for a typical case as follows.\\
\textbf{Case 2}: $p_1$ is odd and $p_2$ is even
We observe the first three equations in the system (3.1) and obtain the following conditions \\
$r_0=0,\frac{1}{2} \qquad$ (i)\\
Since $p_1r_1=k_1 \in \mathbb{Z}$, then $r_1=\frac{m_1}{p_1}, - \lfloor \frac{p_1}{2} \rfloor \leq m_1 \leq \lfloor \frac{p_1}{2} \rfloor,~m_1 \in \mathbb{Z} \qquad$ (ii)\\ As $p_1$ is odd, we have $\lfloor \frac{p_1}{2} \rfloor<\frac{p_1}{2}$, and consequently $r_1$ never equal to $\pm\frac{1}{2}$.\\ 
From $p_2r_2=k_2 \in \mathbb{Z}$, we have $r_2=\frac{m_2}{p_2}, -\lfloor \frac{p_2}{2} \rfloor \leq m_2 \leq \lfloor \frac{p_2}{2} \rfloor,~m_2 \in \mathbb{Z} \qquad$ (iii)\\In particular, as $p_2$ is even, then $\lfloor \frac{p_2}{2} \rfloor=\frac{p_2}{2}$. As a consequence $r_2=\pm\frac{1}{2}$ by choosing $m_2=\pm \lfloor \frac{p_2}{2} \rfloor$. Now, we observe two subcases based on $r_0$.\\

\noindent \textbf{Subcase} $r_0=0$,\\
From $p_2(r_1+r_0)\in \mathbb{Z}$ (The last equation in 3.1), we have $p_2r_1 \in \mathbb{Z}$ then $r_1=\frac{m_2}{p_2}, m_2=-\lfloor \frac{p_2}{2} \rfloor,\dots \lfloor \frac{p_2}{2} \rfloor$.
Therefore, based on the condition (ii) we have 
\begin{align*}
    \{ r_1 \}&=\left\{ \frac{m_1}{p_1}~\bigg\rvert~ m_1=-\left \lfloor \frac{p_1}{2} \right \rfloor,\dots \left\lfloor \frac{p_1}{2} \right\rfloor \right\} \cap \left\{ \frac{m_2}{p_2} ~ \bigg\rvert~ m_2=-\left\lfloor \frac{p_2}{2} \right\rfloor,\dots \left\lfloor \frac{p_2}{2} \right\rfloor \right\}\\
    &=\left\{ \frac{m}{gcd(p_1,p_2)} ~\bigg\rvert~ m=-\left\lfloor \frac{gcd(p_1,p_2)}{2} \right\rfloor,\dots \left\lfloor \frac{gcd(p_1,p_2)}{2} \right\rfloor \right\}
\end{align*}
Furthermore $r_2=r_1+r_0=r_1=\frac{m}{gcd(p_1,p_2)}$. Then for $r_0=0$, we have solutions set\\
\begin{equation*}
    (\tau_0,\tau_1,\tau_2) \sim \left( 0,\frac{m}{gcd(p_1,p_2)}, \frac{m}{gcd(p_1,p_2)}  \right).
\end{equation*}
\textbf{Subcase} $r_0=\pm\frac{1}{2}$\\
From $p_2(r_1+r_0)\in \mathbb{Z}$ (The last equation in 3.1), we have $p_2r_1 + p_2 \cdot \frac{1}{2} \in \mathbb{Z}$. Since $p_2$ is even then $r_1=\frac{m_2}{p_2}, m_2=-\frac{p_2}{2},\dots \frac{p_2}{2}$.
Therefore, based on the condition (ii) we have
\begin{align*}
    \{ r_1 \}&=\left\{ \frac{m_1}{p_1}  ~\bigg\rvert~ m_1=-\left\lfloor \frac{p_1}{2} \right\rfloor,\dots \left\lfloor \frac{p_1}{2} \right\rfloor \right\} \cap \left\{ \frac{m_2}{p_2}  ~\bigg\rvert~ m_2=-\frac{p_2}{2},\dots \frac{p_2}{2} \right\}\\
    &=\left\{ \frac{m}{gcd(p_1,p_2)}  ~\bigg\rvert~ m=-\left\lfloor \frac{gcd(p_1,p_2)}{2} \right\rfloor,\dots \left\lfloor \frac{gcd(p_1,p_2)}{2} \right\rfloor \right\}
\end{align*}
Since $p_1$ is odd, then $gcd(p_1,p_2)$ is an odd number. Hence we have $-\frac{1}{2}< \frac{m}{gcd(p_1,p_2)} < \frac{1}{2} $.\\
Moreover, $r_2=r_1+r_0=\frac{m}{gcd(p_1,p_2)}\pm\frac{1}{2}$. Therefore, by considering (iii) the possible $r_2$ are described by
\begin{align*}
    \{r_2\}&=\left\{ \frac{m}{gcd(p_1,p_2)} \pm \frac{1}{2} \right\} \cap \left\{\frac{m_2}{p_2} ~\bigg\rvert~ m_2=-\frac{p_2}{2},\dots \frac{p_2}{2}\right\}\\
    &=\left\{ \pm\frac{1}{2}\right\}, \text{only satisfied by } m=0, \text{correspond to } r_1=0.
\end{align*}
Then for $r_0=\frac{1}{2}$, we have solutions of the form $\sim \left(\pm \frac{1}{2},0, \pm\frac{1}{2}\right)$.
Hence the complete solution is described by
\begin{equation*}
    (\tau_0,\tau_1,\tau_2)  \sim  \left(0,\frac{m}{gcd(p_1,p_2)},\frac{m}{gcd(p_1,p_2)}\right),\left( \pm\frac{1}{2},0,\pm\frac{1}{2} \right)
\end{equation*}
The similar approach can prove the other cases 1, 3, and 4. $\square$\\
\begin{example}\label{Ex:translation 264}
    Let $\Gamma$ be a space group with $\Gamma_0=(2,6,4)$. Then by Theorem \ref{Thm.3.1} the solution to Frobenius congruence is given by\\ $(\tau_0,\tau_1,\tau_2)\sim(0,0,0),(0,\pm\frac{1}{2},\pm\frac{1}{2}),(\pm\frac{1}{2},0,\pm\frac{1}{2}),(\pm\frac{1}{2},\pm\frac{1}{2},0)$.\\ Or equivalently
    $(r_0,r_1,r_2)=(0,0,0),(0,\pm\frac{1}{2},\pm\frac{1}{2}),(\pm\frac{1}{2},0,\pm\frac{1}{2}),(\pm\frac{1}{2},\pm\frac{1}{2},0)$.
    That means,\\ $(\tau_0,\tau_1,\tau_2)=\xi(0,0,0),\xi(0,\pm\frac{1}{2},\pm\frac{1}{2}),\xi(\pm\frac{1}{2},0,\pm\frac{1}{2}),\xi(\pm\frac{1}{2},\pm\frac{1}{2},0)$, where $0< \xi \in \mathbb{R}$. The suitable value of $\xi$ is determined in the optimal radius computation $\rho$ step.
\end{example}
\subsubsection{Determination of the locally optimal ball packing configuration and its densty for each simply transitive case}
\begin{figure}
    \centering
    \includegraphics[scale=0.33]{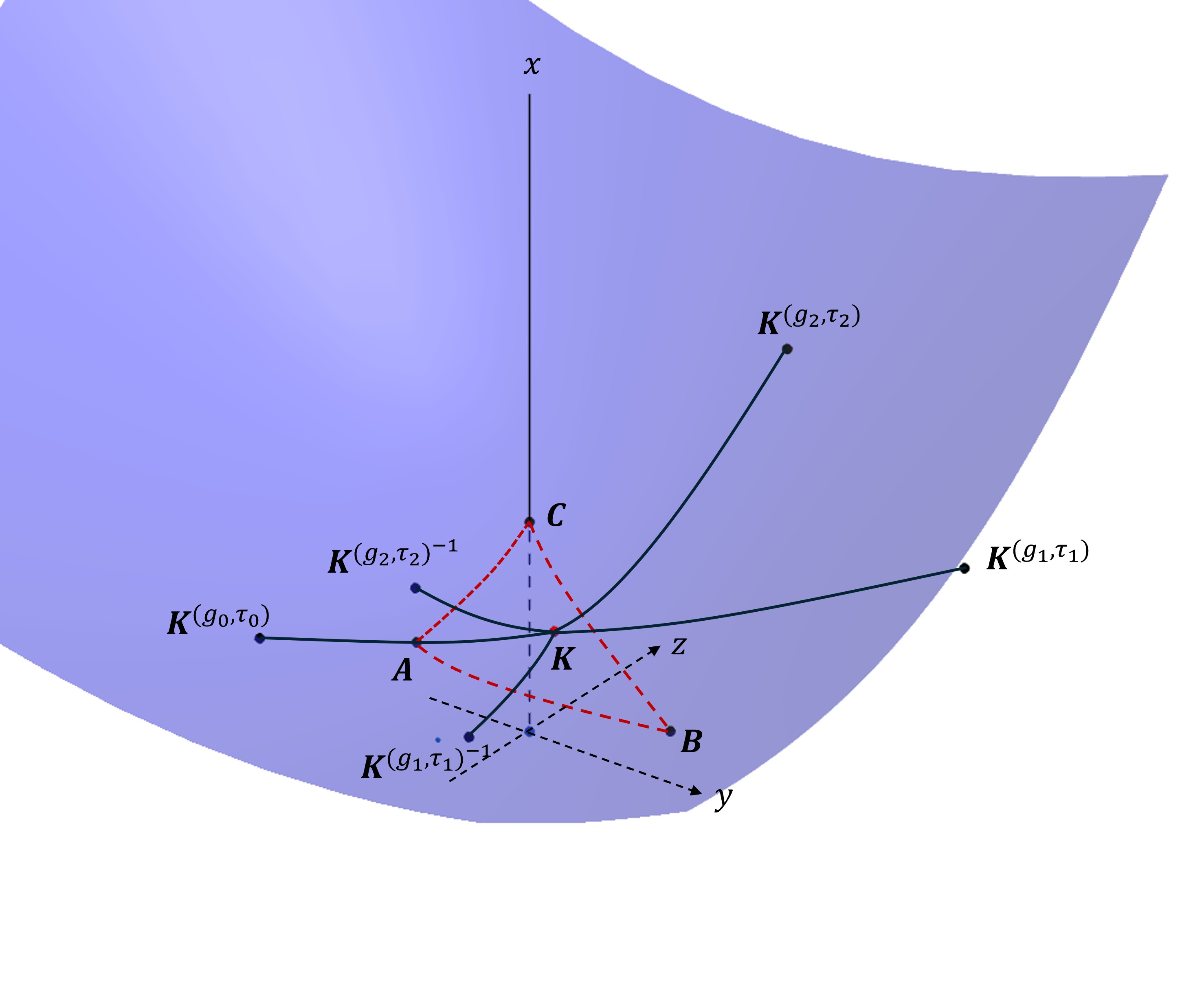}
    \caption{The sketch of rotational point center $A$, $B$, $C$, kernel point $\boldsymbol{K}$ lie on the plane of 0 real fibre coordinate and some neighboring orbit 
    $\boldsymbol{K}^{(g_0,\tau_0)}$, $\boldsymbol{K}^{(g_1,\tau_1)}$, $\boldsymbol{K}^{(g_2,\tau_2)}$}
    \label{Sketch_1}
\end{figure}

We have a complete solution set of Frobenius congruence. Now, we consider the radius of the optimal geodesic ball for each possible parameters $(2,p_1,p_2)$ 
and the corresponding translations $(r_0,r_1,r_2)$ obtained from Theorem \ref{Thm.3.1}, $r_i=\frac{\tau_i}{\xi
}$, for some $\xi>0$. Let $\boldsymbol{K}:=(K,t)$, where $K \in \mathbf{H}^2$ and $t \in \mathbf{R}$, (we may assume that $t=0$) be a kernel point of optimal radius 
$\rho^{opt}$. Then, we consider the neighbouring orbit of $\boldsymbol{K}$ under screw motion group action, i.e $\boldsymbol{K}^{(g_0,\tau_0)}$, $\boldsymbol{K}^{(g_1,\tau_1)}$, 
$\boldsymbol{K}^{(g_2,\tau_2)}$, $\boldsymbol{K}^{(g_0,\tau_0)^{-1}}$, $\boldsymbol{K}^{(g_1,\tau_1)^{-1}}$, $\boldsymbol{K}^{(g_2,\tau_2)^{-1}}$, $\boldsymbol{K}^{(g_0,\tau_0)^{p_0}} 
\equiv \boldsymbol{K}^{(g_1,\tau_1)^{p_1}} \equiv \boldsymbol{K}^{(g_2,\tau_2)^{p_2}} \equiv \boldsymbol{K}^{(\bold{I},\pm 2\rho^{opt})}$, see in Fig. \ref{Sketch_1}.
The necessary condition for optimal congruent geodesic ball packing is that the distance between any two elements of the orbit belonging to the given point $\boldsymbol{K}$ 
is greater than or equal to twice the spherical radius:
\begin{equation}\label{Eq:Key1}
\begin{aligned}
\displaystyle{
    \frac{1}{2}d(\boldsymbol{K},\boldsymbol{K}^{g_0\tau_0})=\frac{1}{2}d(\boldsymbol{K},\boldsymbol{K}^{g_1\tau_1})=\frac{1}{2}d(\boldsymbol{K},\boldsymbol{K}^{g_2\tau_2})=\rho^{opt}}\\
\displaystyle{
    \frac{1}{2}d(\boldsymbol{K},\boldsymbol{K}^{{g_0\tau_0}^{-1}})=\frac{1}{2}d(\boldsymbol{K},\boldsymbol{K}^{{g_1\tau_1}^{-1}})=\frac{1}{2}d(\boldsymbol{K},\boldsymbol{K}^{{g_2\tau_2}^{-1}})=\rho^{opt}}.
\end{aligned}
\end{equation}

\begin{rem} If $\tau_0 \neq 0$ , then the distance of centre geodesic ball and 'upper' center of the neighbor ball is given by $2\rho^{opt}=d(\boldsymbol{K},\boldsymbol{K}^{{(g_0\tau_0)}^{2}})=d((K,t),(K,t+2\tau_0))=2\tau_0=2r_0\cdot \xi=2\cdot \frac{1}{2}\cdot \xi=\xi$. Hence, we have $\xi=2\rho^{opt}$. Otherwise, if $\tau_0=0$, we can choose other $\tau_i \neq 0$, to conclude the similar fact that $\xi=2\rho^{opt}$.
Hence the optimal radius can be formulated as
\begin{equation}\label{Eq:Key}
    \begin{aligned}
    \displaystyle{
        \frac{1}{2}d(\boldsymbol{K},\boldsymbol{K}^{g_0\tau_0})=\frac{1}{2}d(\boldsymbol{K},\boldsymbol{K}^{g_1\tau_1})=\frac{1}{2}d(\boldsymbol{K},\boldsymbol{K}^{g_2\tau_2})=\rho^{opt}}\\
    \displaystyle{\tau_0=r_0 \xi,~ \tau_1=r_1 \xi,~ \tau_2=r_2 \xi,~ \rho^{opt}=\frac{\xi}{2}}.
    \end{aligned}
\end{equation}
\end{rem}
To analyze the system (\ref{Eq:Key}), lets consider the triangle $ABC$ in $\HXR$ space whose the interior angles $\frac{\pi}{2},~\frac{\pi}{p_1},~\frac{\pi}{p_2}$ at vertex $A$, $B$, and $C$ respectively, where the rotational point group will be centered at vertices of triangle $ABC$. Precisely, we choose these vertices point on hyperbolic base plane at real fibre coordinate $t=0$, i.e. $A=(A_0,0)$, $B=(B_0,0)$, $C=(C_0,0)$, where $A_0$, $B_0$, $C_0$ are described in Beltrami-Cayley-Klein model of hyperbolic base plane as follows
\begin{equation}
    \begin{aligned}
        A_0=\left(1,0,-\sqrt{1-\frac{\sin^2{\frac{\pi}{p_2}}}{\cos^2{\frac{\pi}{p_1}}}}\right),~~~~C_0=\left( 1,0,0\right)\\
        B_0=\left(1,\tan{\frac{\pi}{p_2}}\sqrt{1-\frac{\sin^2{\frac{\pi}{p_2}}}{\cos^2{\frac{\pi}{p_1}}}},-\sqrt{1-\frac{\sin^2{\frac{\pi}{p_2}}}{\cos^2{\frac{\pi}{p_1}}}} \right).
    \end{aligned}
\end{equation}
Alternatively, in our projective model as described in Section 2, one could represent the vertices $A,~B,~C$ as
\begin{equation}
\begin{aligned}
    A&=\left[ \left(1,~ \frac{\cos{\frac{\pi}{p_1}}}{\sin{\frac{\pi}{p_2}}},~0,~-\frac{\cos{\frac{\pi}{p_1}}}{\sin{\frac{\pi}{p_2}}}\sqrt{1-\frac{\sin^2\frac{\pi}{p_2}}{\cos^2\frac{\pi}{p_1}}} \right)\right]\\
    B&=\left[ \left(1,~\frac{1}{\tan{\frac{\pi}{p_1}}\tan{\frac{\pi}{p_2}}},~\frac{1}{\tan{\frac{\pi}{p_1}}}\sqrt{1-\frac{\sin^2\frac{\pi}{p_2}}{\cos^2\frac{\pi}{p_1}}},~ -\frac{1}{\tan{\frac{\pi}{p_1}}\tan{\frac{\pi}{p_2}}}\sqrt{1-\frac{\sin^2\frac{\pi}{p_2}}{\cos^2\frac{\pi}{p_1}}} \right)\right]\\
    C&=\left[ \left(1,~1,~0,~0 \right)\right]
\end{aligned}
\end{equation}
Now, let $\boldsymbol{K}$ be the kernel point which in the projective model is described by the polar coordinate: $\boldsymbol{K}=\left[ \left(1,~~\cosh{r},~~\sinh{r}\cos{\alpha},~~\sinh{r}\sin{\alpha} \right) \right]$, $r>0$, $-\pi < \alpha \leq \pi$. Equivalently, $\boldsymbol{K}=(K,0)$, where $K=(1,~\tanh{r}\cos{\alpha},~ \tanh{r}\sin{\alpha})$ in Beltrami-Cayley-Klein model of hyperbolic base plane.\\
By acting screw motion group, we obtain the neighboring orbit of $\boldsymbol{K}$:\\ $\boldsymbol{K}^{(g_0,\tau_0)},~ \boldsymbol{K}^{(g_1,\tau_1)},~ \boldsymbol{K}^{(g_2,\tau_2)}$. In the projective model, these points can be represented as
\begin{equation}
\boldsymbol{K}^{(g_i,\tau_i)}=\boldsymbol{K}^{(g_i, r_i \xi)}:=\left[ \left(1,~\boldsymbol{K}^{(g_i, r_i \xi)}_x,~\boldsymbol{K}^{(g_i, r_i \xi)}_y,~\boldsymbol{K}^{(g_i, r_i \xi)}_z  \right) \right],~~i=0,1,2
\end{equation}
where each entries $\boldsymbol{K}^{(g_i, r_i \xi)}_x,~\boldsymbol{K}^{(g_i, r_i \xi)}_y,~\boldsymbol{K}^{(g_i, r_i \xi)}_z$ are depended on variables $r,~\alpha,~\xi$. Due to the size of the mathematical expressions, we provide these entries in our appendix. In short, we denote $\boldsymbol{K}^{(g_i,\tau_i)}:=\boldsymbol{K}^{(g_i, r_i \xi)}(r,\alpha,\xi)$, for $i=0,1,2$.
\begin{figure}
    \centering
    \includegraphics[scale=0.3]{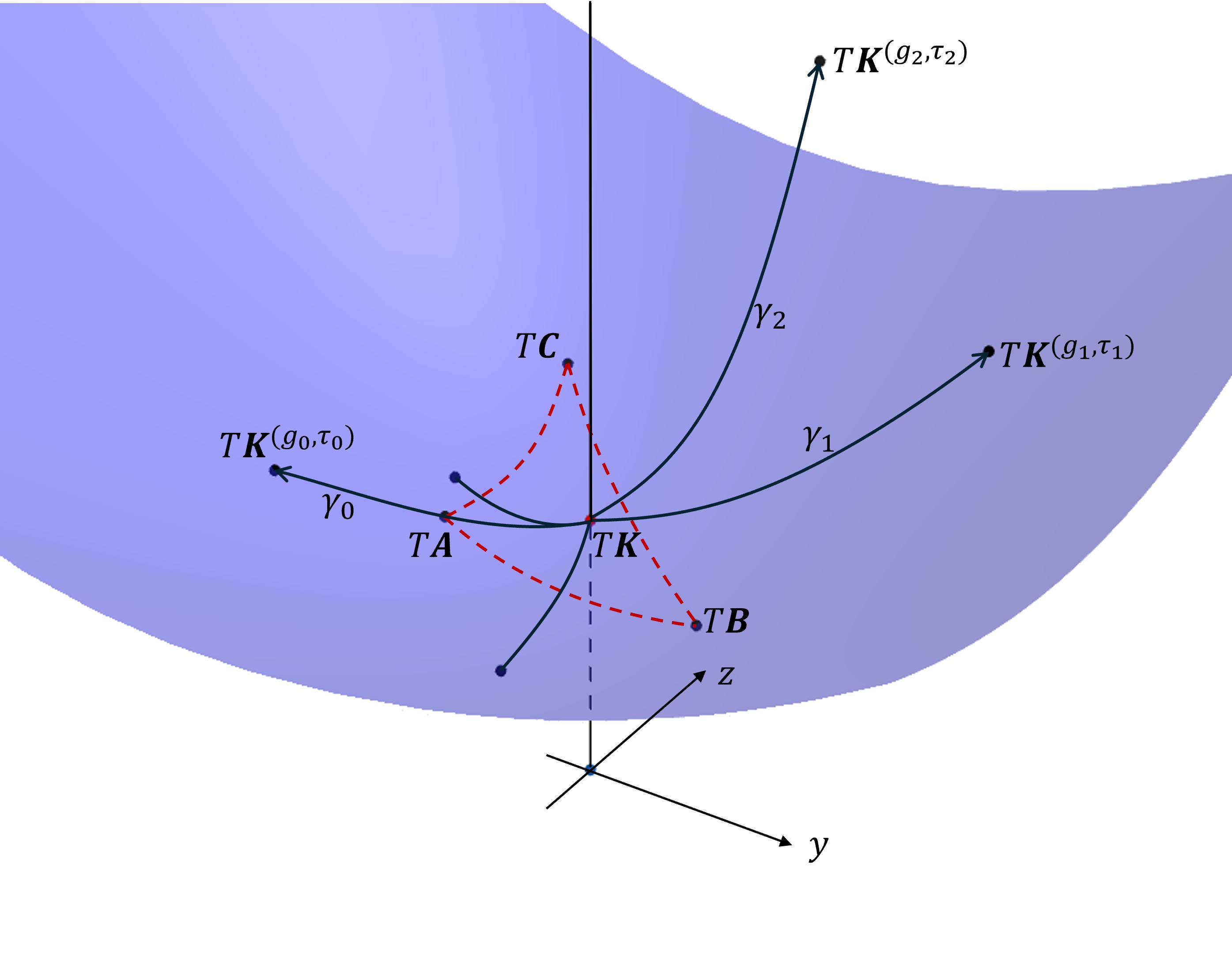}
    \caption{The sketch of the images of rotational centres $A,~B,~C$, kernel $\boldsymbol{K}$ and its nearby orbit $\boldsymbol{K}^{(g_i,\tau_i)}$ under hyperbolic isometry $T$. 
    In this scene, $T\boldsymbol{K}$ lies on the center of model $(1,1,0,0)$. The geodesics $\gamma_i$ emanate from $T\boldsymbol{K}$ to $T\boldsymbol{K}^{(g_i,\tau_i)}$.}
    \label{sketch2}
\end{figure}
To compute the distances $d\left(\boldsymbol{K},\boldsymbol{K}^{(g_i,\tau_i)}\right)$, we consider the geodesic joining $\boldsymbol{K}$ and $\boldsymbol{K}^{(g_i,\tau_i)}$, and measure its length. On the other hand, the geodesic curve starting from the center of the model, $(1,1,0,0)$, with geographical direction $u,v$ is formulated in (\ref{eq:geodesic}). Therefore, it is convenient to place the kernel $\boldsymbol{K}$ to the center of the model $(1,1,0,0)$, by applying an isometry $T$ on $\HXR$ i.e. hyperbolic base translation that takes $\boldsymbol{K}$ to $(1,1,0,0)$, then considering geodesic curve $\gamma_i$ that emanates from $T\boldsymbol{K}=(1,1,0,0)$ to the endpoint $T\boldsymbol{K}^{(g_i,\tau_i)}$, $i=0,1,2$. \\

Hence, we have translated points as follows:
$T\boldsymbol{K}=(1,1,0,0)$, $T\boldsymbol{K}^{(g_i,\tau_i)}$, $i=0,1,2$ (see Fig.\ref{sketch2}). The geodesic starting from $T\boldsymbol{K}=(1,1,0,0)$ to $T\boldsymbol{K}^{(g_i,\tau_i)}$ is formulated as in (\ref{eq:geodesic}), i.e.
{\footnotesize{\begin{equation}
    \gamma_i(\xi)=\left[\left(1,~e^{\xi \cdot \sin{v_i}}\cosh{(\xi \cdot v_i)},~e^{\xi \cdot \sin{v_i}}\sinh{(\xi \cdot v_i)}\cos{u_i},~ e^{\xi \cdot \sin{v_i}}\sinh{(\xi \cdot v_i)}\sin{u_i}\right)\right]
\end{equation}}}
Note that $\gamma_i$ is an arc-length parametrized curve, so the length of $\gamma_i$ as the parameter running from $0$ to $\xi$ is exactly given by $\xi$. 
For some appropriate unique values of $u_i$, $v_i$, and $\xi$, the geodesic curve $\gamma_i$ will reach the point $T\boldsymbol{K}^{(g_i,\tau_i)}$ perfectly.\\
It means that the following equations must hold.
\begin{equation}\label{Eq:4.8}
    \begin{aligned}
        e^{\xi \cdot \sin{v_i}}\cosh{(\xi \cdot v_i)}-T\boldsymbol{K}^{(g_i, r_i \xi)}_x=0\\
        e^{\xi \cdot \sin{v_i}}\sinh{(\xi \cdot v_i)}\cos{u_i}-T\boldsymbol{K}^{(g_i, r_i \xi)}_y=0\\
        e^{\xi \cdot \sin{v_i}}\sinh{(\xi \cdot v_i)}\sin{u_i}-T\boldsymbol{K}^{(g_i, r_i \xi)}_z=0,\\ \text{for } i=0,1,2
    \end{aligned}
\end{equation}
The solution of the system (\ref{Eq:4.8}) is of the form $(u_0,~v_0,~u_1,$ $~v_1,~u_2,~v_2,~r,~\alpha,~\xi)$. 
By unifying all of the equations, we consider the following mapping
\begin{equation}\label{eq:4.9}
\begin{aligned}
    F&:U \subset \mathbf{R}^9 \longrightarrow \mathbf{R}^9,\\
    &:(u_0,~v_0,~u_1,~v_1,~u_2,~v_2,~r,~\alpha,~\xi) \longmapsto F(u_0,~v_0,~u_1,~v_1,~u_2,~v_2,~r,~\alpha,~\xi)
\end{aligned}
\end{equation}
where $U:=\left((-\pi,\pi] \times \left[-\frac{\pi}{2},\frac{\pi}{2}\right]\right)^3 \times (0,\infty) \times (-\pi,\pi] \times (0,\infty)$, and \\
$F:=(f^0_x,~f^0_y,~f^0_z,~f^1_x,~f^1_y,~f^1_z,~f^2_x,~f^2_y,~f^2_z)$ with
\begin{align*}
    f^i_x&:=e^{\xi \cdot \sin{v_i}}\cosh{(\xi \cdot v_i)}-T\boldsymbol{K}^{(g_i, r_i \xi)}_x\\
    f^i_y&:=e^{\xi \cdot \sin{v_i}}\sinh{(\xi \cdot v_i)}\cos{u_i}-T\boldsymbol{K}^{(g_i, r_i \xi)}_y\\
    f^i_z&:=e^{\xi \cdot \sin{v_i}}\sinh{(\xi \cdot v_i)}\sin{u_i}-T\boldsymbol{K}^{(g_i, r_i \xi)}_z,\\ &\text{for } i=0,1,2
\end{align*}
Therefore, we would like to find an appropriate solution of $F=0$ in this setting. Alternatively, $F=0$ could be solved numerically by some methods e.g. by applying the multidimensional Newton-Rhapson scheme. Once it is solved, we get $\xi$, $r$, $\alpha$ then we directly obtain $\rho^{opt}=\frac{\xi}{2}$ and the coordinate of kernel point $\boldsymbol{K}:=(K,0)$, where $K=(1,~\tanh{r}\cos{\alpha},~ \tanh{r}\sin{\alpha})$ in Beltrami-Cayley-Klein model of hyperbolic base plane.
\begin{example}\label{example:264}
    Consider point group $\Gamma_0$ correspond to $(2,6,4)$, and the corresponding triangle $ABC$, see Fig.\ref{Fig.Ex.264} , where $A=\left[\left(1, \frac{\sqrt{6}}{2},0,-\frac{\sqrt{2}}{2} \right) \right]$, $B=\left[\left(1, \sqrt{3}, 1, -1 \right) \right]$, $C=\left[\left(1, 1, 0, 0 \right) \right]$ (on projective model of $\HXR$) are the point rotational center of $\frac{\pi}{2}$, $\frac{\pi}{6}$, $\frac{\pi}{4}$, respectively.
    \begin{figure}
        \centering
        \includegraphics[scale=0.4]{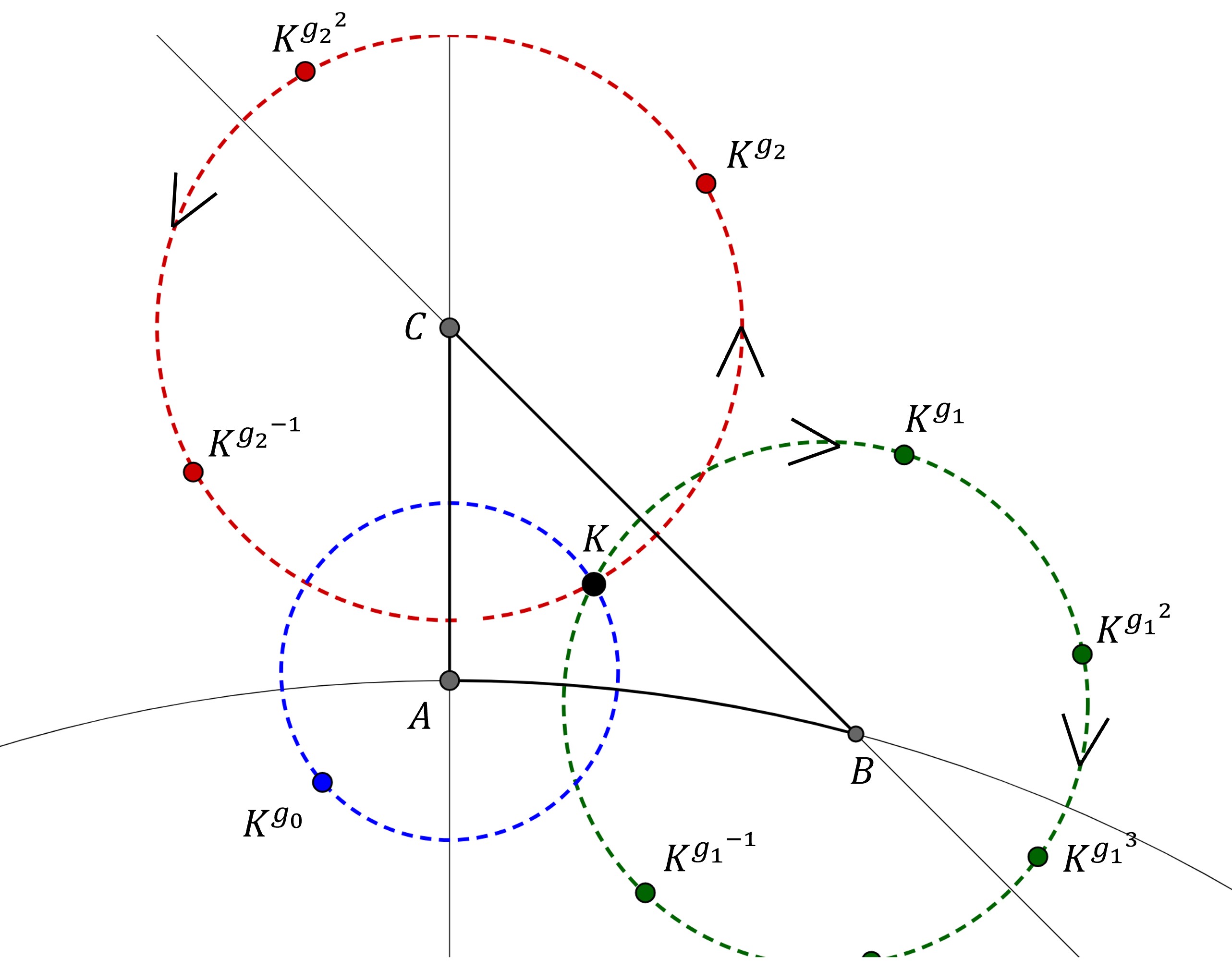}
        \caption{The rotational point group $\Gamma_0$ of $(2,6,4)$ acting on $\bH^2$ (sketched on Poinc\'{a}re disk model) and the orbit of kernel point $\boldsymbol{K}$ under the action $\Gamma_0$.}
        \label{Fig.Ex.264}
    \end{figure}
    We have already seen in example \ref{Ex:translation 264} that the possible translation is\\ $(\tau_0,\tau_1,\tau_2)\sim(0,0,0),(0,\pm\frac{1}{2},\pm\frac{1}{2}),(\pm\frac{1}{2},0,\pm\frac{1}{2}),(\pm\frac{1}{2},\pm\frac{1}{2},0)$.\\ Let's pick one translation part i.e $\sim\left(0,\frac{1}{2},\frac{1}{2}\right)$. The zeros of $F=0$ related to $\sim\left(0,\frac{1}{2},\frac{1}{2}\right)$, is  $\{u_0=-2.46941\dots$, $v_0=0$, $u_1=0.31798\dots$, $v_1=0.523598\dots$, $u_2=1.475468\dots$, $v_2=0.523598$, $r=0.489201\dots$, $\alpha=-0.938184\dots$, $\xi=0.814141\dots\}$. Therefore, the optimum radius of the inscribed geodesic ball is $\rho^{opt}=\frac{\xi}{2}=0.40707\dots$. Moreover from the values $r=0.489201\dots$ and $\alpha=-0.938184\dots$, we obtain the kernel point (center of the geodesic ball) lies on the interior of hyperbolic triangle $ABC$ i.e.
    \begin{equation*}
    \begin{aligned}
        &\boldsymbol{K}=(K,0),\quad \text{where } K=(1,~0.268182,~ -0.365808) \quad \text{(in B-C-K model)}\\
        &\text{or }\boldsymbol{K}=\left[\left( 1,~1.122064,~0.300917,~-0.410460 \right)\right] \quad \text{in the projective model}.
    \end{aligned}
    \end{equation*}
    Furthermore, the volume of geodesic ball $\mathcal{B}_{\boldsymbol{K}}(\rho^{opt})$  and Dirichlet-Voronoi cell $\mathcal{D}(\boldsymbol{K})$ can be directly computed by formula (\ref{Eq:volume geodesic ball}) and (\ref{Eq: Volume of Cell}), that is $\mathrm{Vol}(\mathcal{B}_{\boldsymbol{K}}(\rho^{opt}))=0.28568\dots$ and $\mathrm{Vol}(\mathcal{D}(\boldsymbol{K}))=0.42628\dots$. Finally, we obtain the density packing $\delta=0.67018\dots$.\\
    We perform the same procedure to the other possible translation part $(\tau_0,\tau_1,\tau_2)$. Then, we obtain the optimal radius of their inscribed geodesic ball to the Dirichlet-Voronoi cells formed. The results of our computations are summarized in table \ref{table: simply 264}
    \begin{table}[h!]
    \centering
    \begin{footnotesize}
        \begin{tabular}{||c|c|c|c||}
        \hline
        $(\tau_0,\tau_1,\tau_2)~$&$\rho$ & $\mathrm{vol}(\mathcal{B})$ & $\delta$ \\
        \hline\hline
              $(0,0,0)$ &  $0.35877\dots$ & $0.19510\dots$& $0.51930\dots$\\
        \hline
             $(0,\pm \frac{1}{2},\pm\frac{1}{2})$ & $\bf{0.40707\dots}$ & $\bf{0.28568\dots}$ & $\bf{0.67018\dots}$ \\
        \hline
             $(\pm\frac{1}{2},\pm\frac{1}{2},0)$ & $0.39304\dots$ & $0.25697\dots$ & $0.62433\dots$ \\
        \hline
             $(\pm\frac{1}{2},0,\pm\frac{1}{2})$ & $0.38558\dots$ & $0.24251\dots$ & $0.60060\dots$ \\
        \hline
        \end{tabular}
    \end{footnotesize}
    \caption{The radius, geodesic ball volume and packing density for each possible translation $\tau_0,\tau_1,\tau_2$ in case $(2,6,4)$. 
    The local maximum packing density is $0.67018\dots$.}
    \label{table: simply 264}
\end{table}
\end{example}
\begin{Remark}
Four classes translation part of the form $(\tau_0,\tau_1,\tau_2)\sim(0,0,0)$,\\$(0,\pm\frac{1}{2},\pm\frac{1}{2}),(\pm\frac{1}{2},0,\pm\frac{1}{2}),(\pm\frac{1}{2},\pm\frac{1}{2},0)$ yields equivariant packing. Hence, the packing configurations and their densities are the same, as shown in table 1.
\end{Remark}
\subsubsection{Global maximum packing density for simply transitive cases}
\begin{figure}[h!]
    \centering
    \includegraphics[width=\textwidth]{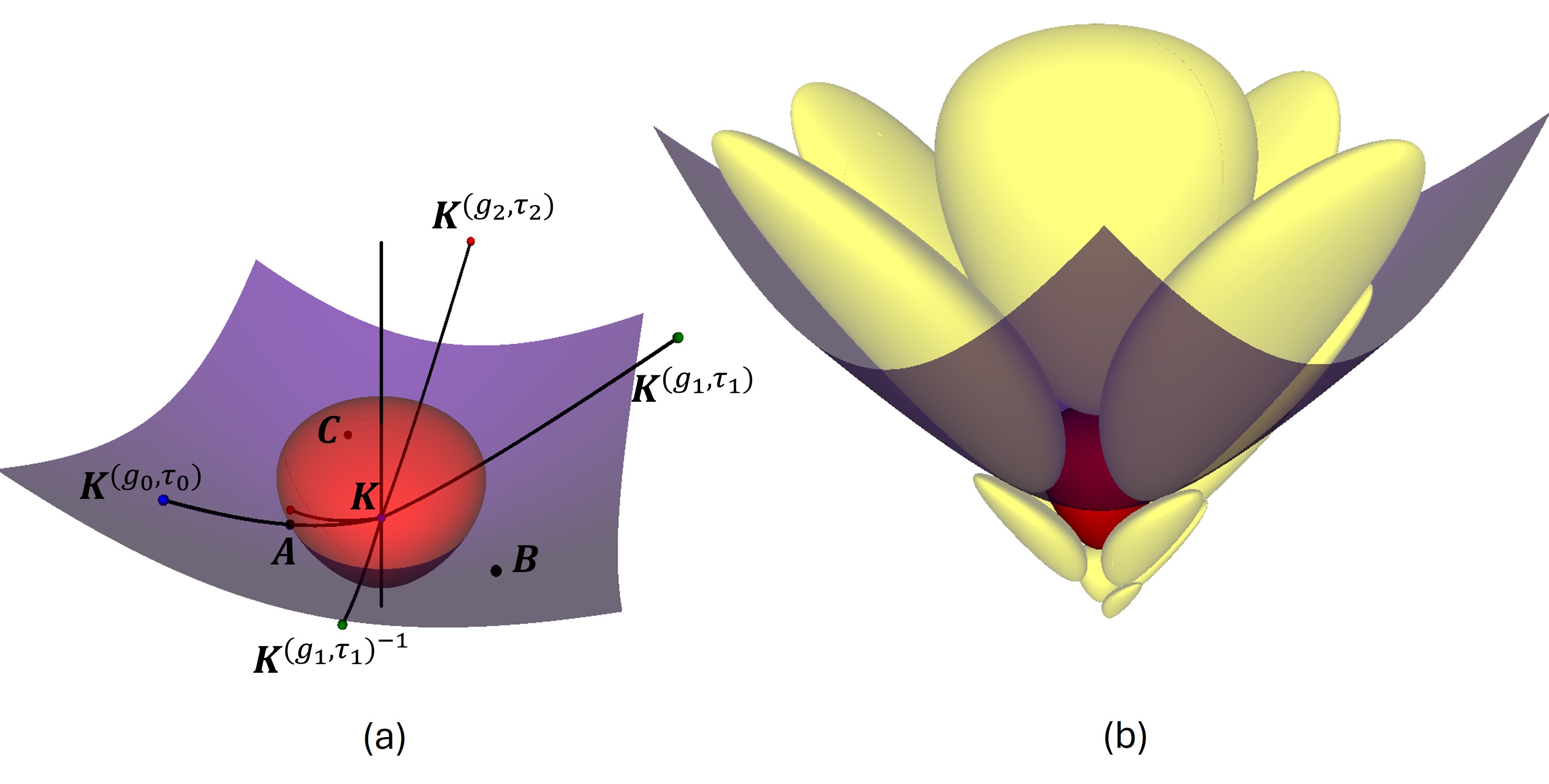}
    \caption{Computer visualization: Optimal packing configurations of the simply-transitive case that attained by $(2,6,4)$ with density $0.67018\dots$, (see Table \ref{tab:simply-transitive 2_p_4}). In this scene, the kernel $\boldsymbol{K}$ lies on $(1,1,0,0)$ and we do not distinguish the translated points by isometry $T$.}
\label{fig:optimum_packing_simply_transitive}
\end{figure}

For each rotational point group $\Gamma_0$ of type $(2,p_1,p_2)$ there exist suitable translation parts $\tau_0,\tau_1,\tau_2$ that forms screw motions group. As shown in example \ref{example:264}, each translation part gives the packing arrangements. Hence, we have the local maximum packing density among them.\\
Now, we consider all types $(2,p_1,p_2)$, and change the parameters $p_1,~p_2$ and all possible translations parts of each type. That is we are looking for the 'global maximum' packing density. By applying the same procedure we obtain all results and present them in tables \ref{tab:simply-transitive 2_p_3}, \ref{tab:simply-transitive 2_p_4}, \ref{tab:simply-transitive 2_p_5} 
for point group generators with given by parameters $(2,p_1,3),~(2,p_1,4),~(2,p_1,5)$, respectively. We do not provide the observations results of $(2,p_1,p_2)$ for the next $p_2 > 5$ here, in fact, their densities are lower than those in $(2,6,4)$ with $\delta^{opt}\approx 0.67018$.\\
Finally, we give a nice visualization of the optimum packing arrangement in Fig.\ref{fig:optimum_packing_simply_transitive}. (b) where it is based on the construction of the optimum kernel obtained $\boldsymbol{K}$ and its nearby orbit under $\Gamma$, i.e. see Fig.\ref{fig:optimum_packing_simply_transitive}. (a).

The proof of the the monotonicity behavior of densities is possible with analytical tools (in another case, see \cite{YSz23-1}). A detailed description of this will be 
discussed in a forthcoming work. Using the above results we obtain the following 
\begin{Theorem}
The optimal packing density configuration of geodesic ball packings generated by screw motions in simply transitive cases is realized 
with parameters $(2,6,4)$, where the optimal density is $\delta^{opt}=0.67018\dots$.
\end{Theorem}

\begin{table}[h!]
    \centering
    \begin{footnotesize}
        \begin{tabular}{||c|c|c|c|c||}
        \hline
        $(p_0,p_1,p_2)$ & $(\tau_0,\tau_1,\tau_2)$&$\rho$ & $\mathrm{vol}(\mathcal{B})$ & $\delta$ \\
        \hline\hline
             $(2,7,3)$   & $(0,0,0)$ &  $0.18773\dots$ & $0.02777\dots$& $0.49454\dots$\\
        \hline
             $(2,8,3)$   &$(\frac{1}{2},\frac{1}{2},0)$ & $0.27040\dots$ & $0.08322\dots$ & $0.58780\dots$ \\
        \hline
             $(2,9,3)$   &$(0,\frac{1}{3},\frac{1}{3})$ & $0.29013\dots$ & $0.10288\dots$ & $0.50792\dots$ \\
        \hline
             $(2,10,3)$   &$(\frac{1}{2},\frac{1}{2},0)$ & $0.33097\dots$ & $0.15298\dots$ & $0.55173\dots$ \\
        \hline
             $(2,11,3)$   &$(0,0,0)$ & $0.31211\dots$ & $0.12819\dots$ & $0.43142\dots$ \\
        \hline
             $(2,12,3)$   &$(\frac{1}{2},\frac{1}{2},0)$ & $0.36067\dots$ & $0.19824\dots$ & $0.52487\dots$ \\
        \hline
            $\vdots$& $\vdots$ &  $\vdots$ & $\vdots$ &$\vdots$ \\
        \hline
            $(2,20,3)$   &$(\frac{1}{2},\frac{1}{2},0)$ & $0.40122\dots$ & $0.27347\dots$ & $0.46491\dots$ \\
         \hline
            $\vdots$& $\vdots$ &  $\vdots$ & $\vdots$ &$\vdots$ \\
        \hline
            $(2,p_1 \rightarrow \infty,3)$   & $(\frac{1}{2},\frac{1}{2},0)$ &$0.42298\dots$ & $0.32081\dots$ & $0.36213\dots$\\
        \hline
        \end{tabular}
    \end{footnotesize}
    \caption{Simply transitive case $(2,p_1,3)$}
    \label{tab:simply-transitive 2_p_3}
\end{table}
\begin{table}[h!]
    \centering
    \begin{footnotesize}
        \begin{tabular}{||c|c|c|c|c||}
        \hline
        $(p_0,p_1,p_2)$ & $(\tau_0,\tau_1,\tau_2)$&$\rho$ & $\mathrm{vol}(\mathcal{B})$ & $\delta$ \\
        \hline\hline
             $(2,5,4)$   & $(\frac{1}{2},0,\frac{1}{2})$ &  $0.30618\dots$ & $0.12099\dots$& $0.62892\dots$ \\
        \hline
             $\textbf{(2,~6,~4)}$   &$(0,\frac{1}{2},\frac{1}{2})$ & $\bf{0.40707\dots}$ & $\bf{0.28568\dots}$ & $\textbf{0.67018\dots}$ \\
        \hline
             $(2,7,4)$   &$(\frac{1}{2},0,\frac{1}{2})$ & $0.42812\dots$ & $0.33273\dots$ & $0.57723\dots$ \\
        \hline
             $(2,8,4)$   &$(0,\frac{1}{2},\frac{1}{2})$ & $0.48221\dots$ & $0.47702\dots$ & $0.62976\dots$ \\
        \hline
             $(2,9,4)$   &$(\frac{1}{2},0,\frac{1}{2})$ & $0.47174\dots$ & $0.44632\dots$ & $0.54208\dots$ \\
        \hline
             $(2,10,4)$   &$(0,\frac{1}{2},\frac{1}{2})$ & $0.51509\dots$ & $0.58265\dots$ & $0.60010\dots$ \\
        \hline
             $(2,11,4)$   &$(\frac{1}{2},0,\frac{1}{2})$ & $0.49294\dots$ & $0.50992\dots$ & $0.51743\dots$ \\
        \hline
            $(2,12,4)$   &$(0,\frac{1}{2},\frac{1}{2})$ & $0.53259\dots$ & $0.64488\dots$ & $0.57813\dots$ \\
        \hline
            $\vdots$& $\vdots$ &  $\vdots$ & $\vdots$ &$\vdots$ \\
        \hline
            $(2,20,4)$   &$(0,\frac{1}{2},\frac{1}{2})$ & $0.55770\dots$ & $0.74181\dots$ & $0.52923\dots$ \\
         \hline
            $\vdots$& $\vdots$ &  $\vdots$ & $\vdots$ &$\vdots$ \\
        \hline
            $(2,p_1 \rightarrow \infty,4)$   & $(0,\frac{1}{2},\frac{1}{2})$ &$0.57167\dots$ & $0.79981\dots$ & $0.44533\dots$\\
        \hline
        \end{tabular}
    \end{footnotesize}
    \caption{Simply transitive case $(2,p_1,4)$}
    \label{tab:simply-transitive 2_p_4}
\end{table}
\begin{table}[h!]
    \centering
    \begin{footnotesize}
        \begin{tabular}{||c|c|c|c|c||}
        \hline
        $(p_0,p_1,p_2)$ & $(\tau_0,\tau_1,\tau_2)$&$\rho$ & $\mathrm{vol}(\mathcal{B})$ & $\delta$ \\
        \hline\hline
             $(2,5,5)$   & $(0,\frac{2}{5},\frac{2}{5})$ &  $0.42561\dots$ & $0.32688\dots$& $0.61116\dots$ \\
        \hline
             $(2,6,5)$   &$(\frac{1}{2},\frac{1}{2},0)$ & $0.48745\dots$ & $0.49292\dots$ & $0.60351\dots$ \\
        \hline
             $(2,7,5)$   &$(0,0,0)$ & $0.47862\dots$ & $0.46632\dots$ & $0.49339\dots$ \\
        \hline
             $(2,8,5)$   &$(\frac{1}{2},\frac{1}{2},0)$ & $0.54450\dots$ & $0.68969\dots$ & $0.57598\dots$ \\
        \hline
             $(2,9,5)$   &$(0,0,0)$ & $0.51171\dots$ & $0.57114\dots$ & $0.47022\dots$ \\
         \hline
             $(2,10,5)$   &$(\frac{1}{2},\frac{1}{2},0)$ & $0.57054\dots$ & $0.79499\dots$ & $0.55441\dots$ \\
        \hline
             $(2,11,5)$   &$(0,0,0)$ & $0.52826\dots$ & $0.62909\dots$ & $0.45323\dots$ \\
        \hline
            $(2,12,5)$   &$(\frac{1}{2},\frac{1}{2},0)$ & $0.58464\dots$ & $0.85633\dots$ & $0.53796\dots$ \\
        \hline
            $\vdots$& $\vdots$ &  $\vdots$ & $\vdots$ &$\vdots$ \\
        \hline
            $(2,20,5)$   &$(\frac{1}{2},\frac{1}{2},0)$ & $0.60513\dots$ & $0.95108\dots$ & $0.50028\dots$ \\
         \hline
            $\vdots$& $\vdots$ &  $\vdots$ & $\vdots$ &$\vdots$ \\
        \hline
            $(2,p_1 \rightarrow \infty,5)$   & $(\frac{1}{2},\frac{1}{2},0)$ &$0.61665\dots$ & $1.00740\dots$ & $0.43334\dots$\\
        \hline
        \end{tabular}
    \end{footnotesize}
    \caption{Simply transitive case $(2,p_1,5)$}
    \label{tab:simply-transitive 2_p_5}
\end{table}
\newpage
\subsection{Multiply transitive cases}
Analogously with the simply transitive case, we investigate the kernel 
$\boldsymbol{K}$ in $\HXR$ and the orbit $\boldsymbol{K}^{\Gamma}$ such that the 
density of corresponding ball packing $\mathcal{B}^{\Gamma}(\boldsymbol{K})$ is 
maximal. But, now we consider that $\Gamma_{\boldsymbol{K}} \neq \bold{I}$, where 
it happens if we choose $\boldsymbol{K}$ on the rotational axis through point 
$A,~B$ or $C$. For instance, if we put $\boldsymbol{K}$ on the rotational 
axis through the point $A$, then the rotation $(\bg_0,0)$ belongs to the stabilizer 
$\Gamma_{\boldsymbol{K}}$. The shape of the fundamental domain will change 
significantly as well as the corresponding possible translation part. 

We provide all possible translations parts in the next subsection (Theorem \ref{thm_frobenius_multiply_case}).
\subsubsection{Possible translation part}
\begin{thm}[Frobenius congruences in multiply transitive cases]\label{thm_frobenius_multiply_case}
    For each $(2,p_1,p_2)$ the solutions to Frobenius congruence is satisfied by the following
    
    \begin{enumerate}
        \item For $\boldsymbol{K}=A$, $(\tau_0,\tau_1,\tau_2)\sim (0,\frac{m}{gcd(p_1,p_2)},\frac{m}{gcd(p_1,p_2)})$,\\ where $-\left\lfloor \frac{gcd(p_1,p_2)}{2} \right\rfloor \leq m \leq \left\lfloor \frac{gcd(p_1,p_2)}{2} \right\rfloor,~m \in \mathbb{Z}$.\\
        \item For $\boldsymbol{K}=B$\\
        $(\tau_0,\tau_1,\tau_2)\sim \begin{cases}
            (0,0,0), \left( \pm \frac{1}{2},0,\pm \frac{1}{2} \right),~ \text{if}~p_2~ \text{even}.\\
            (0,0,0),~ \text{if}~p_2~ \text{odd}.
        \end{cases}$
        \item For $\boldsymbol{K}=C$\\
        $(\tau_0,\tau_1,\tau_2)\sim \begin{cases}
            (0,0,0), \left( \pm \frac{1}{2},\pm \frac{1}{2},0 \right),~ \text{if}~p_1~ \text{even}.\\
            (0,0,0),~ \text{if}~p_1~ \text{odd}.
        \end{cases}$
    \end{enumerate}
\end{thm}
\textbf{Proof:}\\
Applying similar observation as in the proof of Theorem \ref{Thm.3.1} we obtain the proof of this theorem analogously. 
$\square$
\subsubsection{Determination of the locally optimal ball packing configuration and
its density for each simply transitive case}
Similarly to the approach used in simply transitive cases, we move the kernel 
point $\boldsymbol{K}$ to $(1,1,0,0)$ by simply applying a 
hyperbolic isometry $T$, since we assumed that $\boldsymbol{K}$ has the same real 
fibre coordinate as the center of the model. 
Then, we consider the geodesics starting from the translated kernel 
(now, it is the starting point of the geodesics (see (2.8)) to the nearby 
transformed orbit points of $T\boldsymbol{K}^{\Gamma}$. 
By solving the systems (see Eq.\ref{Eq:4.8}) resulting from this consideration, 
we completely obtain all the 
information needed, i.e. $\rho^{opt}$, $\boldsymbol{K}^{opt}$, 
and the geodesics $\gamma_i$.\\
\begin{figure}[h!]
    \centering
    \includegraphics[scale=0.35]{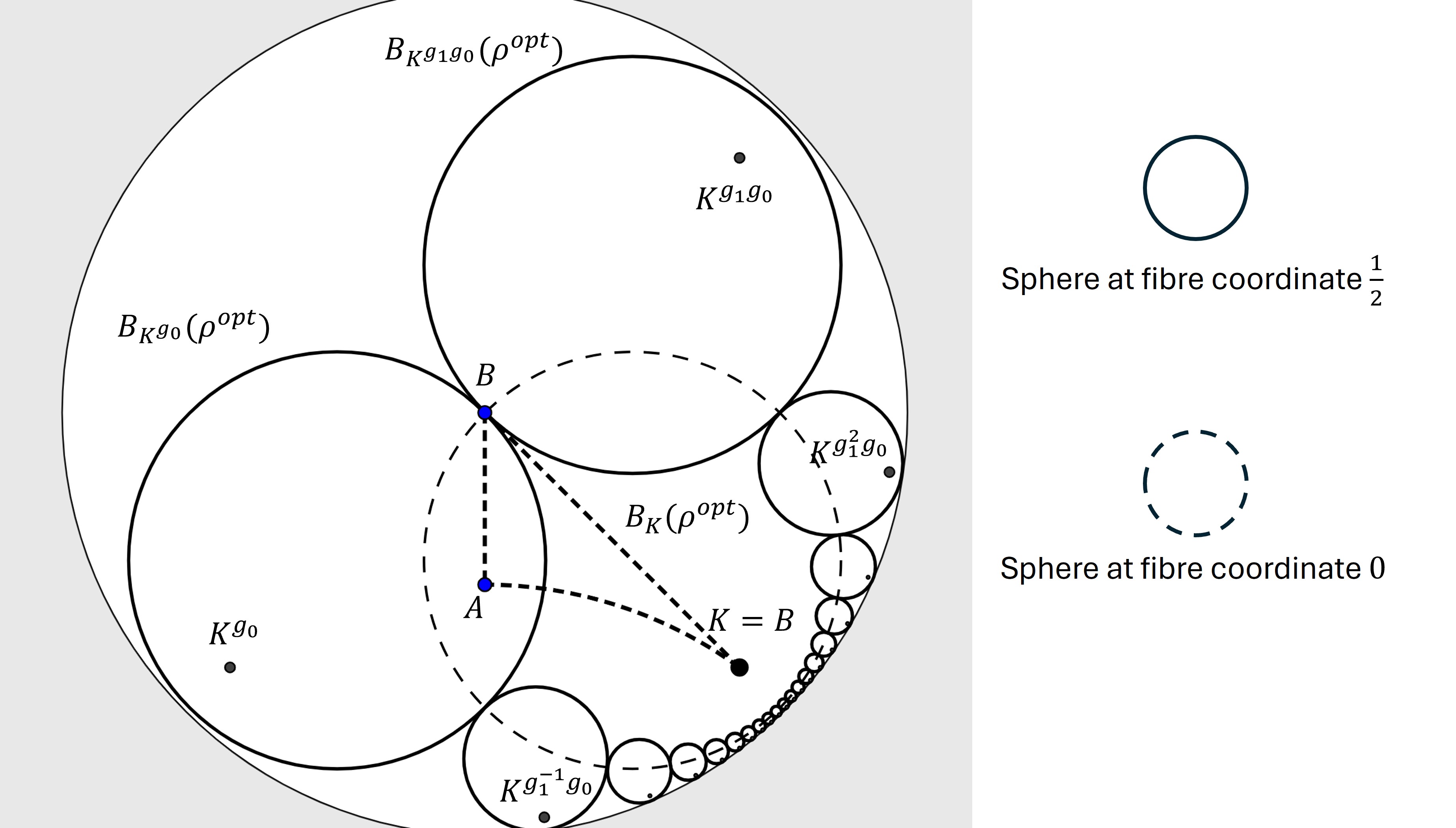}
    \caption{The projection onto Poinc\'{a}re base model $\bold{H}^2$: The sphere centered at $\boldsymbol{K}=B$ of radius $\rho^{opt}$ (with dashed curve), and its nearby images under $\Gamma$ of generator related to $(2,20,4)$ (indicated with full curve). }
    \label{fig:enter-label}
\end{figure}
\subsubsection{The global maximum packing densities for multiply transitive cases}
In general, for every point group generator $\Gamma_0$ of the parameter 
$(2, p_1, p_2)$, there exist three different Kernel choices, each corresponding to a 
rotational center: $A,~B$, and $C$. Furthermore, for every kernel point chosen 
$\boldsymbol{K}$, there exist appropriate translations part that admits lattice in 
the real fibre, as stated in Theorem \ref{thm_frobenius_multiply_case}. 
Hence we formulate the global maximum as in (\ref{Eq:Global_maximum}), 
but the kernel $\boldsymbol{K}$ should be chosen from vertices $A, B, C$, 
not freely from the interior of the triangle.\\

In particular, for some generators $\Gamma_0$ of $(2,p_1,p_2)$, the optimal 
densities are strictly monotonically increasing, but for some finite terms sequences 
of certain chosen kernel points. e.g. for $(2,p_1,4)$, $(2,p_1,6)$, see Tables 6 and 8. 
It shows that by choosing kernel point at $B$, 
then we get the densities as monotonically increasing finite sequences.\\

Precisely, this limitation is due to the condition that $2\rho^{opt}> d(\boldsymbol{K}^{g_1^n,g_0},\boldsymbol{K}^{g_1^{n+1},g_0})$. That means the geodesic balls $\displaystyle B_{\boldsymbol{K}^{(g_1^n,g_0)}}(\rho^{opt})$ and $\displaystyle B_{\boldsymbol{K}^{(g_1^{n+1},g_0)}}(\rho^{opt})$ are overlapping. Hence, we could not take these cases as a ball packing. For instance, on the point group of parameter $(2,p_1,4)$: we have strictly increasing optimum density sequence, for $p_1=5,\dots 20$. But, if we consider $p_1 > 20$, we have overlapping geodesic balls. Therefore, for $p_1>20$, we should choose other kernel points, instead of $B$. Consequently, the global optimum packing density is attained by the group whose parameter $(2,20,4)$, with the optimum density $\approx 0.80529$, see Table 6 and Figure 5. A similar situation occurs in $(2,p_1,6)$, where the monotonicity is satisfied when $p_1\leq 62$, see Table 8. However, the optimum density is less than those in $(2,p_1,4)$. In general, the optimum packing density generated by the further groups of parameter $(2,p_1,p_2)$, where $p_2>6$, $p_2$ even number, are less than $\approx 0.80529$. We do not include complete computation results on the table.\\

The proof of the monotonicity behavior of densities is possible with analytical 
tools (in another case, see \cite{YSz23-1}). A detailed description of this will be 
discussed in a forthcoming work. Using the above results we obtain the following 
\begin{Theorem}
The optimal packing density configuration of congruent geodesic ball packings generated by screw motions in multiply transitive cases is realized 
with parameters $(2,20,4)$, where the optimal density is $\delta^{opt}=0.80529\dots$.
\end{Theorem}
%
%
\begin{table}[h!]\label{tab:multiply_2p3}
    \centering
    \begin{footnotesize}
        \begin{tabular}{||c|c|c|c|c|c||}
        \hline
        $(p_0,p_1,p_2)$ &Optimum &Optimum& $\rho^{opt}$ & $\mathrm{vol}(\mathcal{B}(K))$  &$\delta^{opt}$ \\
        $~$ &Kernel &Translation& ~ &  & \\
        \hline\hline
            $(2,7,3)$   & $B$ &$(0,0,0)$& $0.54527\dots$ & $0.69267\dots$ & $0.60653\dots$ \\
        \hline
            $(2,8,3)$   & $B$ &$(0,0,0)$& $0.76428\dots$ & $1.94411\dots$ & $0.607262\dots$ \\
        \hline
             $(2,9,3)$   & $B$ &$(0,0,0)$& $0.92753\dots$ & $3.53909\dots$& $\bold{0.607267\dots}$ \\
        \hline
             $(2,10,3)$   & $B$&$(0,0,0)$&$1.06127\dots$ & $5.39521\dots$  &$0.606823\dots$ \\
        \hline
             $(2,11,3)$   &$B$&$(0,0,0)$& $1.17585\dots$ & $7.46309\dots$  & $0.60608\dots$ \\
        \hline
             $(2,12,3)$   &$B$&$(0,0,0)$& $1.27668\dots$ & $9.70891\dots$  & $0.60516\dots$ \\
        \hline
            $\vdots$&$\vdots$ &$\vdots$ & $\vdots$ &  $\vdots$  & $\vdots$ \\
        \hline
            $(2,20,3)$   &$B$&$(0,0,0)$ & $1.82969\dots$ & $31.96254\dots$  &$0.59576\dots$\\
        \hline
            $\vdots$&$\vdots$ &$\vdots$ & $\vdots$ &  $\vdots$  & $\vdots$ \\
        \hline
            $(2,p_1 \rightarrow \infty,3)$   &$C$&$(\pm\frac{1}{2},\pm\frac{1}{2},0)$ & $0.63428\dots$ & $1.09790\dots$  &$0.27548\dots$\\
        \hline
        \end{tabular}
        \caption{Multi transitive cases $(2,p_1,3)$: The maximum density is attained at $(2,9,3)$ with density $\delta^{opt}=0.607267\dots$, where the kernel point is $B$ and the translation part $(0,0,0)$}
    \end{footnotesize}
\end{table}
\begin{table}[h!]\label{tab:multiply_2p4}
    \centering
    \begin{footnotesize}
        \begin{tabular}{||c|c|c|c|c|c||}
       \hline
        $(p_0,p_1,p_2)$ &Optimum &Optimum& $\rho^{opt}$ & $\mathrm{vol}(\mathcal{B}(K))$  &$\delta^{opt}$ \\
        $~$ &Kernel &Translation& ~ &  & \\
        \hline\hline
            $(2,5,4)$   & $B$ &$(\pm \frac{1}{2},0,\pm\frac{1}{2})$& $0.72384\dots$ & $1.64498\dots$ & $0.72337\dots$ \\
        \hline
            $(2,6,4)$   & $B$ &$(\pm\frac{1}{2},0,\pm\frac{1}{2})$& $1.01772\dots$ & $4.72953\dots$ & $0.73962\dots$ \\
        \hline
             $(2,7,4)$   & $B$ &$(\pm\frac{1}{2},0,\pm\frac{1}{2})$& $1.23599\dots$ & $8.75091\dots$& $0.75121\dots$ \\
        \hline
             $(2,8,4)$   & $B$&$(\pm \frac{1}{2},0,\pm\frac{1}{2})$&$1.41361\dots$ & $13.50190\dots$  &$0.76007\dots$ \\
        \hline
             $(2,9,4)$   &$B$&$(\pm \frac{1}{2},0,\pm\frac{1}{2})$& $1.56467\dots$ & $18.85517\dots$  & $0.76716\dots$ \\
        \hline
             $(2,10,4)$   &$B$&$(\pm \frac{1}{2},0,\pm\frac{1}{2})$& $1.69666\dots$ & $24.72283\dots$  & $0.77303\dots$ \\
        \hline
             $(2,11,4)$   &$B$&$(\pm \frac{1}{2},0,\pm\frac{1}{2})$&$1.81413\dots$ & $31.04013\dots$ & $0.77804\dots$\\
        \hline
            $(2,12,4)$   &$B$&$(\pm \frac{1}{2},0,\pm\frac{1}{2})$&$1.92013\dots$ & $37.75731\dots$ & $0.78240\dots$ \\
        \hline
            $\vdots$&$\vdots$ &$\vdots$ & $\vdots$ &  $\vdots$  & $\vdots$ \\
        \hline
            $(2,18,4)$   &$B$&$(\pm \frac{1}{2},0,\pm\frac{1}{2})$ & $2.40393\dots$ & $84.67149\dots$  &$0.80082\dots$\\
        \hline
            $(2,19,4)$   &$B$&$(\pm \frac{1}{2},0,\pm\frac{1}{2})$ & $2.46761\dots$ & $93.39040\dots$  &$0.80312\dots$\\
        \hline
            $\bold{(2,20,4)}$   &$B$&$(\pm \frac{1}{2},0,\pm\frac{1}{2})$ & $\bold{2.52789\dots}$ & $\bold{102.32545\dots}$  &$\bold{0.80529\dots}$\\
        \hline
            $(2,21,4)$   &$B$&$(0,0,0)$ & $2.18922\dots$ & $60.06928\dots$  &$0.54587\dots$\\
        \hline
            $(2,22,4)$   &$B$&$(0,0,0)$ & $2.28600\dots$ & $70.30425\dots$  &$0.54385\dots$\\
        \hline
             $(2,23,4)$   &$B$&$(0,0,0)$ & $2.33104\dots$ & $75.53280\dots$  &$0.54285\dots$\\
        \hline
            $(2,24,4)$   &$B$&$(0,0,0)$ & $2.37412\dots$ & $80.82932\dots$  &$0.54185\dots$\\
        \hline
            $\vdots$&$\vdots$ &$\vdots$ & $\vdots$ &  $\vdots$  & $\vdots$ \\
        \hline
            $(2,p_1 \rightarrow \infty,4),p_1 \text{even} $   &$A$&$(0,\pm \frac{1}{2},\pm \frac{1}{2})$ & $0.76034\dots$ & $1.91344\dots$  &$0.40051\dots$\\
            $(2,p_1 \rightarrow \infty,4),p_1 \text{odd} $   &$A$&$(0,0,0)$ & $0.65847\dots$ & $1.23095\dots$  &$0.29752\dots$\\
        \hline
        \end{tabular}
        \caption{Multi transitive cases $(2,p_1,4)$: The maximum density is attained at $(2,20,4)$ with density $\delta^{opt}=0.80529\dots$, where the kernel point is $B$ and the translation part $(\pm \frac{1}{2},0,\pm\frac{1}{2})$}
    \end{footnotesize}
\end{table}

\begin{table}[h!]
    \centering
    \begin{footnotesize}
        \begin{tabular}{||c|c|c|c|c|c||}
        \hline
        $(p_0,p_1,p_2)$ &Optimum &Optimum& $\rho^{opt}$ & $\mathrm{vol}(\mathcal{B}(K))$  &$\delta^{opt}$ \\
        $~$ &Kernel &Translation& ~ &  & \\
        \hline\hline
            $(2,5,5)$   & $A$ &$(0,\pm\frac{2}{5},\pm\frac{2}{5})$& $0.57897\dots$ & $0.83129\dots$ & $0.57128\dots$ \\
        \hline
            $(2,6,5)$   & $C$ &$(\pm\frac{1}{2},\pm \frac{1}{2},0)$& $1.08335\dots$ & $5.75699\dots$ & $\bold{0.63431\dots}$ \\
        \hline
             $(2,7,5)$   & $B$ &$(0,0,0)$& $1.23499\dots$ & $8.72815\dots$& $0.51127\dots$ \\
        \hline
             $(2,8,5)$   & $C$&$(\pm\frac{1}{2},\pm \frac{1}{2},0)$&$1.18247\dots$ & $7.59768\dots$  &$0.58434\dots$ \\
        \hline
             $(2,9,5)$   &$B$&$(0,0,0)$& $1.50608\dots$ & $16.61921\dots$  & $0.51653\dots$ \\
        \hline
             $(2,10,5)$   &$C$&$(\pm\frac{1}{2},\pm \frac{1}{2},0)$& $1.22545\dots$ & $8.51436\dots$  & $0.55289\dots$ \\
         \hline
             $(2,11,5)$   &$B$&$(0,0,0)$& $1.71621\dots$ & $25.69831\dots$  & $0.51807\dots$ \\
        \hline
             $(2,12,5)$   &$C$&$(\pm\frac{1}{2},\pm \frac{1}{2},0)$& $1.24813\dots$ & $9.02911\dots$  & $0.53138\dots$ \\
        \hline
            $\vdots$&$\vdots$ &$\vdots$ & $\vdots$ &  $\vdots$  & $\vdots$ \\
        \hline
            $(2,20,5)$   &$B$&$(0,0,0)$ & $2.32684\dots$ & $75.03233\dots$  &$0.51321\dots$\\
        \hline
            $\vdots$&$\vdots$ &$\vdots$ & $\vdots$ &  $\vdots$  & $\vdots$ \\
        \hline
            $(2,p_1 \rightarrow \infty,5)$   &$C$&$(\pm\frac{1}{2},\pm\frac{1}{2},0)$ & $1.29808\dots$ & $10.24237\dots$  &$0.41859\dots$\\
        \hline
        \end{tabular}
        \caption{Multi transitive cases $(2,p_1,5)$: The maximum density is attained at $(2,6,5)$ with density $\delta^{opt}=0.63431\dots$, where the kernel point is $C$ and the translation part $(\pm\frac{1}{2},\pm\frac{1}{2},0)$}
    \end{footnotesize}
    \label{tab:multiply_2p5}
\end{table}
\begin{table}[h!]
    \centering
    \begin{footnotesize}
        \begin{tabular}{||c|c|c|c|c|c||}
        \hline
        $(p_0,p_1,p_2)$ &Optimum &Optimum& $\rho^{opt}$ & $\mathrm{vol}(\mathcal{B}(K))$  &$\delta^{opt}$ \\
        $~$ &Kernel &Translation& ~ &  & \\
        \hline\hline
            $(2,6,6)$   & $B$ &$(\pm \frac{1}{2},0,\pm\frac{1}{2})$& $1.32353\dots$ & $10.90425\dots$ & $0.65561\dots$ \\
        \hline
            $(2,7,6)$   & $B$ &$(\pm\frac{1}{2},0,\pm\frac{1}{2})$& $1.51801\dots$ & $17.05731\dots$ & $0.67063\dots$ \\
        \hline
             $(2,8,6)$   & $B$ &$(\pm\frac{1}{2},0,\pm\frac{1}{2})$& $1.68241\dots$ & $24.03027\dots$& $0.68197\dots$ \\
        \hline
             $(2,9,6)$   & $B$&$(\pm \frac{1}{2},0,\pm\frac{1}{2})$&$1.82523\dots$ & $31.69581\dots$  &$0.69094\dots$ \\
        \hline
             $(2,10,6)$   &$B$&$(\pm \frac{1}{2},0,\pm\frac{1}{2})$& $1.95167\dots$ & $39.96131\dots$  & $0.69830\dots$ \\
        \hline
             $(2,11,6)$   &$B$&$(\pm \frac{1}{2},0,\pm\frac{1}{2})$& $2.06521\dots$ & $48.75648\dots$  & $0.70451\dots$ \\
        \hline
             $(2,12,6)$   &$B$&$(\pm \frac{1}{2},0,\pm\frac{1}{2})$&$2.16831\dots$ & $58.02620\dots$ & $0.70985\dots$\\
        \hline
            $\vdots$&$\vdots$ &$\vdots$ & $\vdots$ &  $\vdots$  & $\vdots$ \\
        \hline
            $(2,60,6)$   &$B$&$(\pm \frac{1}{2},0,\pm\frac{1}{2})$ & $4.03967\dots$ & $759.64262\dots$  &$0.78758\dots$\\
        \hline
            $(2,61,6)$   &$B$&$(\pm \frac{1}{2},0,\pm\frac{1}{2})$ & $4.05878\dots$ & $777.39291\dots$  &$0.78836\dots$\\
        \hline
            $\bold{(2,62,6)}$   &$B$&$(\pm \frac{1}{2},0,\pm\frac{1}{2})$ & $\bold{4.07757\dots}$ & $\bold{795.23140\dots}$  &$\bold{0.78913\dots}$\\
        \hline
            $(2,63,6)$   &$B$&$(0,0,0)$ & $3.54729\dots$ & $413.04460\dots$  &$0.46329\dots$\\
        \hline
            $(2,64,6)$   &$C$&$(\pm \frac{1}{2},\pm \frac{1}{2},0)$ & $1.51908\dots$ & $17.09714\dots$  &$0.46984\dots$\\
        \hline
            $(2,65,6)$   &$B$&$(0,0,0)$ & $3.57857\dots$ & $429.73373\dots$  &$0.46239\dots$\\
        \hline
            $(2,66,6)$   &$C$&$(\pm \frac{1}{2},\pm \frac{1}{2},0)$ & $1.51918\dots$ & $17.10070\dots$  &$0.46921\dots$\\
        \hline
            $\vdots$&$\vdots$ &$\vdots$ & $\vdots$ &  $\vdots$  & $\vdots$ \\
        \hline
            $(2,p_1 \rightarrow \infty,6), p_1 \text{even}$   &$C$&$(\pm \frac{1}{2},\pm \frac{1}{2},0)$ & $1.52069\dots$ & $17.15693\dots$  &$0.44890\dots$\\
            $(2,p_1 \rightarrow \infty,6), p_1 \text{odd}$   &$C$&$(\pm \frac{1}{2},\pm \frac{1}{2},0)$ & $1.31695\dots$ & $10.73024\dots$  &$0.32418\dots$\\
        \hline
        \end{tabular}
        \caption{Multi transitive cases $(2,p_1,6)$: The maximum density is attained at $(2,62,6)$ with density $\delta^{opt}={0.78913\dots}$,  with kernel point $B$ and the translation part $(\pm \frac{1}{2},0,\pm\frac{1}{2})$}
    \end{footnotesize}
    \label{tab:multiply_2p6}
\end{table}
\newpage
\subsection{Concluding Remarks}
We have just demonstrated that the geodesic ball packing generated by the screw motion group yields numerous configurations and achieves a very high packing density. It would be interesting to investigate the geodesic ball packing constructed by other point group generators. Several further questions merit consideration: What are the possible translation components related to Frobenius congruence and the group structure itself? Would the optimal density be higher than that observed in the current study? We are optimistic that our method can be applied to these general situations.\\
Thurston geometries with non-constant curvature are full-fledged alternatives to those with constant curvature, and the 
investigation of classical questions is just as interesting (for example, due to material structure aspects) 
as in the case of classical geometries. For all these reasons, their investigation may become the focus of further research.

\section{Appendix}
\begin{footnotesize}
\begin{align*}
\boldsymbol{K}^{(g_0,\tau_0)}_x&=\frac{e^{\xi \cdot r_0}}{\sin^{2}\frac{\pi}{p_2}}\left|2\cos^{2}\frac{\pi}{p_1}\left(\cosh{r}+\sinh{r}\sin{\alpha}\sqrt{1-\frac{\sin^{2}\frac{\pi}{p_2}}{\cos^{2}\frac{\pi}{p_1}}} \right)-\sin^{2}\frac{\pi}{p_2}\cosh{r}\right|\\
    \boldsymbol{K}^{(g_0,\tau_0)}_y&=-e^{\xi \cdot r_0}\sinh{r}\cos{\alpha}\\
    \boldsymbol{K}^{(g_0,\tau_0)}_z&=-\frac{2e^{\xi \cdot r_0}}{\sin^{2}\frac{\pi}{p_2}}\left(\cos^{2}\frac{\pi}{p_1} \cosh{r}\sqrt{1-\frac{\sin^{2}\frac{\pi}{p_2}}{\cos^{2}\frac{\pi}{p_1}}}+\sinh{r}\sin{\alpha}\left(\cos^2{\frac{\pi}{p_1}}-\frac{1}{2}\sin^2{\frac{\pi}{p_2}}\right) \right)\\
\boldsymbol{K}^{(g_1,\tau_1)}_x&=\frac{e^{\xi \cdot r_1}}{\sin^2{\frac{\pi}{p_2}}}\left| 2\cos^{2}\frac{\pi}{p_1}\left(\cosh{r}+\sinh{r}\sin{\alpha}\sqrt{1-\frac{\sin^{2}\frac{\pi}{p_2}}{\cos^{2}\frac{\pi}{p_1}}} \right)-\sin^{2}\frac{\pi}{p_2}\cosh{r} \right|\\
\boldsymbol{K}^{(g_1,\tau_1)}_y&=\frac{2e^{\xi \cdot r_1}}{\sin^2{\frac{\pi}{p_2}}}\Bigg(\cosh{r}\cos{\frac{2\pi}{p_2}}\cos^2{\frac{\pi}{p_1}}\sqrt{1-\frac{\sin^2{\frac{\pi}{p_2}}}{\cos^2{\frac{\pi}{p_1}}}}-\sinh{r}\Big(\sin{\frac{\pi}{p_2}}\cos^2{\frac{\pi}{p_2}}\sin{\big(\frac{\pi}{p_2}-\alpha \big)}\\&-\sin{\alpha}\sin{\frac{2\pi}{p_2}}\cos{\frac{2\pi}{p_2}}-\frac{1}{2}\cos{\alpha}\sin^2{\frac{\pi}{p_2}}\Big)\Bigg)\mathrm{csgn}\Bigg(2\sinh{r}\sin{\alpha}\cos^2{\frac{\pi}{p_1}}\sqrt{1-\frac{\sin^2{\frac{\pi}{p_2}}}{\cos^2{\frac{\pi}{p_1}}}}\\&-\cosh{r}\Big(\sin^2{\frac{\pi}{p_2}}-2\cos^2{\frac{\pi}{p_1}} \Big)\Bigg)\\
\boldsymbol{K}^{(g_1,\tau_1)}_z&=-\frac{2e^{\xi \cdot r_1}}{\sin^2{\frac{\pi}{p_2}}}
\Bigg(\cosh{r}\cos{\frac{2\pi}{p_2}}\cos{\frac{\pi}{p_1}}\sqrt{\cos^2{\frac{\pi}{p_2}}-\sin^2{\frac{\pi}{p_1}}}\\
&-\sinh{r}\Big( \cos^3{\frac{\pi}{p_2}}\sin{\big(\frac{\pi}{p_2}-\alpha \big)}+\cos^2\frac{\pi}{p_2}\sin{\alpha}\big(\frac{3}{2}-2\cos^2{\frac{\pi}{p_1}} \big)-\frac{1}{2}\cos{\alpha}\sin{\frac{2\pi}{p_2}}\\
&+\frac{1}{2}\sin{\alpha}\cos{\frac{2\pi}{p_1}}\Big)\Bigg)\mathrm{csgn}\Bigg(2\sinh{r}\sin{\alpha}\cos^2{\frac{\pi}{p_1}}\sqrt{1-\frac{\sin^2{\frac{\pi}{p_2}}}{\cos^2{\frac{\pi}{p_1}}}}\\
&-\cosh{r}\left(\sin^2{\frac{\pi}{p_2}}-2\cos^2{\frac{\pi}{p_1}} \right) \Bigg)\\
%
\boldsymbol{K}^{(g_2,\tau_2)}_x&=e^{\xi \cdot r_2}\cosh{r}\\
\boldsymbol{K}^{(g_2,\tau_2)}_y&=e^{\xi \cdot r_2}\sinh{r}\cos{\left(\alpha+\frac{2\pi}{p_2}\right)}\\
\boldsymbol{K}^{(g_2,\tau_2)}_z&=e^{\xi \cdot r_2}\sinh{r}\sin{\left(\alpha+\frac{2\pi}{p_2}\right)}
\end{align*}
\end{footnotesize}

\begin{thebibliography}{MPSz98}
%
%
\bibitem{F01}
{Farkas,~Z.~J.}
The classification of $\SXR$ space groups,
\emph{Beitr. Algebra Geom.,}
{\bf42} (2001), 235--250.
%
\bibitem{FTL} Fejes~T\'oth,~L.
Regular Figures,
\textit{Macmillan (New York)}, 1964.
%
\bibitem{G--K--K} Fejes~T\'oth,~L.~---~Fejes~T\'oth,~G.~---~Kuperberg,~W.
Lagerungen: Arrangements in the Plane, on the Sphere, and in Space
{\bf 360} Grundlehren der mathematischen Wissenschaften
{\textit Springer Nature}, (2023), ISBN	3031218000, 9783031218002.
%
\bibitem{M}
  {Macbeath,~A.~M}
  The classification of non-Euclidean plane crystallographic groups. 
  \emph{Can. J. Math.,}
 {\bf19} (1967), 1192--1205.
 %
\bibitem{M97}
{Moln{\'a}r,~E.}
The projective interpretation of the eight 3-di\-men\-sional homogeneous geometries. 
\emph{Beitr. Algebra Geom.,}
{\bf38} No.~2 (1977) 261--288.
%
\bibitem{MSz}
{Moln{\'a}r,~E.~---~Szirmai,~J.}
Symmetries in the 8 homogeneous 3-geometries.
\textit{Symmetry Cult. Sci.,}
{\bf 21/1-3} (2010), 87--117.
%
\bibitem{MSz18} Moln\'ar,~E.~---~Szirmai,~J.
\textit{Top dense hyperbolic ball packings and coverings for complete Coxeter orthoscheme groups},
{Publications de l'Institut Mathématique}, {\bf 103(117)}  (2018), 129--146, DOI: 10.2298/PIM1817129M.
%
\bibitem{MSz12}
{Moln{\'a}r,~E.~---~Szirmai,~J.}
Classification of $\SOL$ lattices.
\textit{Geom. Dedicata,}
{\bf 161/1} (2012), 251-275.
%
\bibitem{MSzV}
Moln{\'a}r,~E.~---~Szirmai,~J.~---~Vesnin,~A.
Projective metric realizations of cone-manifolds with singularities along 2-bridge knots and links.
{\it J. Geom.,}  {\bf 95} (2009), 91-133.
%
\bibitem{MSz14}
Moln{\'a}r,~E.~---~Szirmai,~J.
Volumes and geodesic ball packings to the regular prism tilings in $\SLR$ space.
{\it Publ. Math. Debrecen}  {\bf 84}(1-2) (2014), 189--203.
%
\bibitem{MSzV17}
Moln{\'a}r~E.~---~Szirmai~J.~---~Vesnin~A.,
{Geodesic and Translation Ball Packings
Generated by Prismatic Tesselations of  the Universal Cover of $\SLR$},
{\it Results in Math.}  {\bf 71)} (2017) 623--642.
%
\bibitem{N17}
N{\'e}meth,~L.
Pascal pyramid in the space $\HXR$.  
{\it Mathematical Communications}  {\bf 22} (2017), 211--225.
%
\bibitem{PSSz10}
{Pallagi,~J.~---~Schultz,~B.~---~Szirmai,~J.}
Visualization of geodesic curves, spheres and equidistant surfaces in $\SXR$ space.
\emph{KoG}, {\bf 14} (2010), 35--40.
%
\bibitem{PSSz11-2}
{Pallagi,~J.~---~Schultz,~B.~---~Szirmai,~J.}
Equidistant surfaces in $\HXR$ space.
\emph{KoG}, {\bf 15}, (2011), 3-6.
%
\bibitem{PSz12}
{Pallagi,~J.~---~Szirmai,~J.}
Visualization of the Dirichlet-Voronoi cells in $\SXR$ space.
\emph{Pollack Periodica}, {\bf 7} Supp 1, 95--104 (2012), DOI: 10.1556/Pollack.7.2012.S.9.
%
%
\bibitem{S}
Scott,~P.
The geometries of 3-manifolds. 
{\it Bull. London Math. Soc.} {\bf 15} (1983), 401--487.
%
\bibitem{Sz07}
Szirmai,~J.
The densest geodesic ball packing by a type of $\NIL$ lattices.
{\it Beitr. Algebra Geom.} {\bf 48}(2) (2007), 383--398.
%
\bibitem{Sz13-1}
Szirmai,~J.
A candidate to the densest packing with equal balls in the Thurston geometries. 
{\it Beitr. Algebra Geom.,} {\bf 55}(2) (2014), 441--452.
%
%
\bibitem{Sz11-1}
{Szirmai,~J.}
Simply transitive geodesic ball packings to $\mathbf{S^2\times R}$ space groups generated by glide reflections,
{\emph {Ann. Mat. Pur. Appl.}}, {\bf 193/4} (2014), 1201-1211, DOI: 10.1007/s10231-013-0324-z.
%
\bibitem{Sz11-2}
{Szirmai,~J.}
Geodesic ball packings in $\SXR$ space for generalized Coxeter space groups.
\emph{\it Beitr. Algebra Geom.,}
{\bf 52}, (2011), 413 -- 430.
 %
\bibitem{Sz12-5}
{Szirmai,~J.}
Geodesic ball packings in $\HXR$ space for generalized Coxeter space groups.
\textit{Math. Commun.}, {\bf 17/1} (2012), 151--170.
%
\bibitem{Sz20}
Szirmai,~J.
Interior angle sums of geodesic triangles in $\SXR$ and $\HXR$ geometries.  
{\it Bul. Acad. de Stiinte Republicii Mold. Mat.,} {\bf 93(2)} (2020), 44--61.
%
\bibitem{Sz21}
Szirmai,~J.
Apollonius surfaces, circumscribed spheres of tetrahedra, Menelaus' and Ceva's theorems in $\SXR$ and $\HXR$ geometries. 
{\it Q. J. Math.,} (2021), DOI: 10.1093/qmath/haab038.
%
\bibitem{Sz22}
Szirmai,~J.
On Menelaus' and Ceva's theorem in $\NIL$ geometry. 
{\it Submitted Manuscript,} (2022), arXiv: 2110.08877.
%
\bibitem{Sz22-3}
{Szirmai,~J.:}
Classical Notions and Problems in Thurston Geometries,
\emph{International Electronic Journal of Geometry},
{\bf 16} No.2 (2023), 608--643, doi: 10.36890/IEJG.1221802, arXiv: 2203.05209.
%
\bibitem{Sz23}
{Szirmai,~J.:}
Fibre-like cylinders, their packings and coverings in $\SLR$ space,
\emph{Results Math.}, (2024), doi: 10.1007/s00025-024-02152-0, arXiv: 2306.05721.
%
\bibitem{Sz06-1} Szirmai,~J. 
The $p$-gonal prism tilings and their optimal hypersphere packings in the hyperbolic
3-space,
\textit{Acta Math. Hungar.}, {\bf{111 (1-2)}} (2006), 65--76.
%
\bibitem{Sz06-2} Szirmai,~J. The regular prism tilings and their optimal hyperball packings in the hyperbolic $n$-space,
\textit{Publ. Math. Debrecen}, {\bf{69 (1-2)}} (2006), 195--207.
%
\bibitem{Sz13-3} Szirmai,~J. The optimal hyperball packings related to the smallest compact arithmetic $5$-orbifolds,
\textit{Kragujevac J. Math.} {\bf 40(2)} (2016), 260-270, DOI:10.5937/KgJMath1602260S.
%
\bibitem{Sz13-4} Szirmai,~J. The least dense hyperball covering to the regular prism tilings in the hyperbolic $n$-space,
\textit{Ann. Mat. Pur. Appl.} {\bf 195/1} (2016), 235--248, DOI: 10.1007/s10231-014-0460-0.
%
\bibitem{T}
Thurston,~W.~P. (and Levy,~S. editor),
{\it Three-Dimensional Geometry and Topology}.  Princeton University Press,  Princeton, New Jersey, vol. {\bf 1} (1997).
%
\bibitem{W06}
Weeks,~J.~R.
{Real-time animation in hyperbolic, spherical, and product geometries.}
\textit{A.~Pr\'ekopa and E.~Moln\'ar, (eds.).
Non-Euclidean Geometries, J\'anos Bolyai Memorial Volume,
Mathematics and Its Applications}, Springer (2006) Vol.~{\bf 581}, 287--305.
%
\bibitem{YA2023}
Yahya,~A.
{On Problem of Best Circle to Discontinuous Groups in Hyperbolic Plane.}
\textit{Mathematical Communications}, ~{\bf 28 (1)}, 121--140.
%
\bibitem{YSz23-1}
{Yahya, A.~--~Szirmai, J.:}
Geodesic ball packings generated by rotations and monotonicity behavior of their densities in $\mathbf{H}^2\!\times\!\mathbf{R}$ space.
\emph{Submitted manuscript}, (2023).
%
\end{thebibliography}
\end{document}